\documentclass[ijoc,nonblindrev]{simpleFormat} 
\usepackage{enumerate}

\OneAndAHalfSpacedXII 

\usepackage{natbib}
 \bibpunct[, ]{(}{)}{,}{a}{}{,}%
 \def\bibfont{\small}%
\TheoremsNumberedThrough     

\EquationsNumberedThrough    

\MANUSCRIPTNO{JOC-2019-09-OA-214} 

\newcommand{\E}[1]{{\mathrm{E}\left[#1\right]}}

\newcommand{\PP}[1]{\mathrm{P}\left( #1 \right)}

\newcommand{\T}{{\mathrm{T}}}

\newcommand{\edit}[1]{{#1}}

\begin{document}


\RUNAUTHOR{Daw and Pender}

\RUNTITLE{Matrix Calculations for Moments of Markov Processes}

\TITLE{Matrix Calculations for Moments of Markov Processes}

\ARTICLEAUTHORS{%
\AUTHOR{Andrew Daw}
\AFF{School of Operations Research and Information Engineering, Cornell University, \EMAIL{amd399@cornell.edu}}
\AUTHOR{Jamol Pender}
\AFF{School of Operations Research and Information Engineering, Cornell University, \EMAIL{jjp274@cornell.edu}}
} 

\ABSTRACT{%
Matryoshka dolls, the traditional Russian nesting figurines, are known world-wide for each doll's encapsulation of a sequence of smaller dolls. In this paper, we identify a large class of Markov process whose moments are easy to compute by exploiting the structure of a new sequence of nested matrices we call \textit{Matryoshkhan matrices}. We characterize the salient properties of Matryoshkhan matrices that allow us to compute these moments in closed form at a specific time without computing the entire path of the process.  This speeds up the computation of the Markov process moments significantly in comparison to traditional differential equation methods, which we demonstrate through numerical experiments.  Through our method, we derive explicit expressions for both transient and steady-state moments of this class of Markov processes. We demonstrate the applicability of this method through explicit examples such as shot-noise processes, growth-collapse processes, ephemerally self-exciting processes, and affine stochastic differential equations from the finance literature.  We also show that we can derive explicit expressions for the self-exciting Hawkes process, for which finding closed form moment expressions has been an open problem since its introduction in 1971. \edit{In general, our techniques can be used for any Markov process for which the infinitesimal generator of an arbitrary polynomial is itself a polynomial of equal or lower order.}
}%



\maketitle

%
\section{Introduction.}\label{intro} 

In recently studying the intensity of Markovian Hawkes process, originally defined in \citet{hawkes1971spectra}, we have been interested in computing all the moments of this process. \edit{In surveying the literature for this process, there does not seem to be any closed form transient solutions at the fourth order or higher (see Proposition 5 of \citet{gao2018limit} for moments one through three), and both steady-state solutions and ordinary differential equations have only been available up to the fourth moment, see \citet{da2014Hawkes, errais2010affine}. Similarly, \citet{ait2015modeling} give expressions for the fourth transient moment of a self-exciting jump-diffusion model up to squared error in the length of time, and one could simplify these expressions to represent the Markovian Hawkes intensity with the same error.} The standard methodology for finding moments is to differentiate the moment generating function to obtain the moments, however, this is intractable for practical reasons, see for example \citet{errais2010affine}. The problem of finding the moments of the Hawkes process is also the subject of the recent interesting research in \citet{cui2019elementary,cui2019moments}, works that are concurrent and independent from this one. In \citet{cui2019elementary}, the authors propose a new approach for calculating moments that they construct from elementary probability arguments and also relate to the infinitesimal generator. Like the infinitesimal generator, this new methodology produces differential equations that can be solved algebraically or numerically to yield the process moments, and the authors provide closed form transient expressions up to the second moment. \citet{cui2019moments} extends this methodology to cases of Gamma decay kernels. In other recent previous works, \citet{daw2018queues, koops2018infinite}, the authors have identified the differential equation for an arbitrary moment of the Hawkes process, although the closed form solutions for these equations have remained elusive and prompted closer investigation. Upon inspecting the differential equation for a given moment of the Hawkes process intensity, one can notice that this expression depends on the moments of lower order. Thus, to compute a given moment one must solve a system of differential equations with size equal to the order of the moment, meaning one must at least implicitly solve for all the lower order moments first. This same pattern occurs in \citet{cui2019elementary}. Noticing this nesting pattern leads one to wonder: what other processes have moments that follow this structure?

In this paper, we explore this question by identifying what exactly this nesting structure is. In the sequel, we will define a novel sequence of matrices that captures this pattern. Just as Matryoshka dolls -- the traditional Russian nesting figurines -- stack inside of one another, these matrices are characterized by their encapsulation of their predecessors in the sequence. Hence, we refer to this sequence as \textbf{Matryoshkan matrices}. As we will show, these matrices can be used to describe the linear system of differential equations that arise in solving for the moments of the Hawkes process, as well as the moments of a large class of other Markov processes. In fact, the only assumption we make on these processes is that their moments satisfy differential equations that do not depend on any higher order moments. As we will demonstrate through detailed examples, this includes a wide variety of popular stochastic processes, such as It\^o diffusions and shot noise processes. By utilizing this nesting structure we are able to solve for the moments of these processes in closed form. By comparison to traditional methods of solving these systems of differential equations numerically, the advantage of the approach introduced herein is the fact that the moments can be computed at a specific point in time rather than on a path through time. This yields a methodology that is both efficient and precise.

\edit{
This methodology also has the potential to be quite relevant in practice. Of course, these techniques can be used to efficiently calculate the commonly used first four moments, thus obtaining the mean, variance, skewness, and kurtosis. Moreover though, let us note that the higher moment calculations are also of practical use. For example, these higher moments can be used in Markov-style concentration inequalities, as the higher order should improve the accuracy of the tail bounds. To that end, one can also use the vector of moments to approximate generating functions such as the moment generating function of Laplace-Stieltjes transforms. This could then be used to characterize the stationary distribution of the process, for example, or to provide approximate calculations of quantities such as the cumulative distribution function through transform methods. The calculation of moments can also be highly relevant for many applications in mathematical finance. For example, these techniques may hold great potential for polynomial processes, see for example \citet{filipovic2019polynomial}. One could also expect this efficiently calculated vector of moments to be of use in estimation through method of moment techniques. Again in this case, access to higher order moments should improve fit.
}

\edit{
The remainder of this paper is organized as follows. In Section~\ref{secPrelims}, we introduce Matryoshkan matrix sequences and identify some of their key properties. In Section~\ref{secMoments} we use these matrices to find the moments of a large class of general Markov processes. We also give specific examples. In Section~\ref{secNum}, we demonstrate the numerical performance of this method in comparison to traditional differential equation techniques. In Section~\ref{secConc}, we conclude. Throughout the course of this study, we make the following contributions:
\begin{enumerate}[i)]
\item We define a novel class of matrix sequences that we call {Matryoshkan} matrix sequences for their nesting structure. We identify key properties of these matrices such as their inverse and matrix exponentials.
\item Through these Matryoshkan matrices, we solve for closed form expressions for the moments of a large class of Markov processes. Furthermore, we demonstrate the general applicability of this technique through application to notable stochastic processes including Hawkes processes, shot noise process, It\^o diffusions, growth-collapse processes, and linear birth-death-immigration processes. In the case of the Hawkes process and growth-collapse processes this resolves an open problem, as closed form expressions of these general transient moments were not previously known in the literature.
\item We compare the precision and computation time of our methodology to numerically solving the underlying differential equations. In observing empirical superiority of the Matryoshkan matrix approach, we demonstrate the efficiency of calculating the moments at a given point, rather than on a path through time.
\end{enumerate}
}


\section{Matryoshkan Matrix Sequences}\label{secPrelims}

\edit{
For the sake of clarity,  let us begin this section by  introducing general notation patterns we will use throughout this paper. Because of the heavy use of matrices in this work, we reserve boldface upper case variables for these objects, such as $\mathbf{I}$ for the identity matrix. Similarly we let boldface lower case variables be vectors, such as $\mathbf{v}$ for the vector of all ones or $\mathbf{v}_i$ being the unit vector in the $i^\text{th}$ direction. One can assume that all vectors are column vectors unless otherwise noted. Scalar terms will not be bolded. A special matrix that we will use throughout this work is the diagonal matrix, which we denote $\mathbf{diag}(\mathbf{a})$, which is a square matrix with the values of the vector $\mathbf{a}$ along its diagonal and zeros otherwise. We will also make use of a generalization of this, denoted $\mathbf{diag}(\mathbf{a}, k)$, which instead contains the values of $\mathbf{a}$ on the $k^\text{th}$ off-diagonal, with negative $k$ being below the diagonal and positive $k$ being above.
}

Let us now introduce a sequence of matrices that will be at the heart of this work. We begin as follows: consider a sequence of lower triangular matrices $\{\mathbf{M}_n, n \in \mathbb{Z}^+\}$ such that
\begin{align}\label{matryoshkan}
\mathbf{M}_{n}
=
\begin{bmatrix}
\mathbf{M}_{n-1} & \mathbf{0}_{n-1 \times 1} \\
\mathbf{m}_{n} & m_{n,n}
\end{bmatrix}
,
\end{align}
where $\mathbf{m}_{n} \in \mathbb{R}^{n-1}$ is a row vector, $m_{n,n} \in \mathbb{R}$, and $\mathbf{M}_1 = m_{1,1}$, an initial value. Taking inspiration from Matryoshka dolls, the traditional Russian nesting dolls, we will refer to these objects as \textit{Matryoshkan} matrices. Using their nested and triangular structures, we can make four quick observations of note regarding Matryoshkan matrices.

\begin{proposition}\label{matryoshkanProp}
Each of the following statements is a consequence of the definition of Matryoshkan matrices given by Equation~\ref{matryoshkan}:
\begin{enumerate}[i)]
\item If $\mathbf{X}_n \in \mathbb{R}^{n \times n}$ and $\mathbf{Y}_n \in \mathbb{R}^{n \times n}$ are both Matryoshkan matrix sequences, then so are $\mathbf{X}_n + \mathbf{Y}_n$ and $\mathbf{X}_n\mathbf{Y}_n$.
\item
If $m_{i,i} \ne 0$ for all $i \in \{1, \dots , n\}$ then the Matryoshkan matrix $\mathbf{M}_n\in \mathbb{R}^{n \times n}$ is nonsingular. Moreover, the inverse of $\mathbf{M}_{n}$ is given by the recursion
\begin{align}\label{mminv}
\mathbf{M}_n^{-1}
=
\begin{bmatrix}
\mathbf{M}_{n-1}^{-1} & \mathbf{0}_{n-1 \times 1} \\
-\frac{1}{m_{n,n}}\mathbf{m}_{n}\mathbf{M}_{n-1}^{-1} & \frac{1}{m_{n,n}}
\end{bmatrix}
.
\end{align}
\item
If $m_{i,i} \ne 0$ for all $i \in \{1, \dots , n\}$ and are all distinct then the matrix exponential of the Matryoshkan matrix $\mathbf{M}_n\in \mathbb{R}^{n \times n}$ multiplied by $t \in \mathbb{R}$ follows the recursion
\begin{align}\label{mmexp}
e^{\mathbf{M}_{n}t}
=
\begin{bmatrix}
e^{\mathbf{M}_{n-1}t} & \mathbf{0}_{n-1 \times 1} \\
\mathbf{m}_{n}
\left(\mathbf{M}_{n-1} - m_{n,n}\mathbf{I}\right)^{-1}
\left(e^{\mathbf{M}_{n-1}t} - e^{m_{n,n} t}\mathbf{I}\right)
&
e^{m_{n,n} t}
\end{bmatrix}
.
\end{align}
\item
If $m_{i,i} \ne 0$ for all $i \in \{1, \dots , n\}$ and are all distinct then the matrices $\mathbf{U}_n \in \mathbb{R}^{n \times n}$ and $\mathbf{D}_n \in \mathbb{R}^{n \times n}$ are such that
$$
\mathbf{M}_n \mathbf{U}_n = \mathbf{U}_n \mathbf{D}_n
$$
for the Matryoshkan matrix $\mathbf{M}_n \in \mathbb{R}^{n \times n}$ when defined recursively as
\begin{align}\label{mmeig}
\mathbf{U}_n
=
\begin{bmatrix}
\mathbf{U}_{n-1} & \mathbf{0}_{n-1 \times 1} \\
\mathbf{m}_{n}\left(\mathbf{D}_{n-1} - m_{n,n}\mathbf{I} \right)^{-1} \mathbf{U}_{n-1} & 1
\end{bmatrix}
,\,
\mathbf{D}_n
=
\begin{bmatrix}
\mathbf{D}_{n-1} & \mathbf{0}_{n-1 \times 1} \\
\mathbf{0}_{1 \times n-1} & m_{n,n}
\end{bmatrix}
.
\end{align}
\end{enumerate}
\end{proposition}

One can pause to note that in some sense any lower triangular matrix could be considered Matryoshkan, or at least be able to satisfy these properties. However, we note that some of the most significant insights we can gain from the Matryoshkan structure are the recursive implications available for sequences of matrices. Moreover, it is the combination of the nested relationship of consecutive matrices and the lower triangular structure that enables us to find these patterns. We will now see how this notion of Matryoshkan matrix sequences and the associated properties above can be used to specify element-wise solutions to a sequence of differential equations.

\begin{lemma}\label{mmlemma}
Let $\mathbf{M}_n \in \mathbb{R}^{n \times n}$, $\mathbf{c}_n \in \mathbb{R}^{n}$, and $\mathbf{s}_n(t) : \mathbb{R}^+ \to \mathbb{R}^n$ be such that
$$
\mathbf{M}_n
=
\begin{bmatrix}
\mathbf{M}_{n-1} & \mathbf{0}_{n-1\times 1} \\
\mathbf{m}_{n} & m_{n,n}\\
\end{bmatrix}
,
\quad
\mathbf{c}_{n}
=
\begin{bmatrix}
\mathbf{c}_{n-1}\\
c_n
\end{bmatrix}
,
\quad
\mathrm{ and }
\quad
\mathbf{s}_{n}(t)
=
\begin{bmatrix}
\mathbf{s}_{n-1}(t)\\
s_n(t)
\end{bmatrix}
$$
where $\mathbf{m}_{n} \in \mathbb{R}^{n-1}$ is a row vector, $\mathbf{c}_{n-1} \in \mathbb{R}^{n-1}$, $s_n(t) \in \mathbb{R}$, and $\mathbf{M}_1 = m_{1,1}$. Further, suppose that
$$
\frac{\mathrm{d}}{\mathrm{d}t}\mathbf{s}_{n}(t)
=
\mathbf{M}_{n}\mathbf{s}_{n}(t)
+
\mathbf{c}_{n} .
$$
Then, if $m_{k,k} \ne 0$ for all $k \in \{1, \dots, n\}$, the vector function $\mathbf{s}_n(t)$ is given by
\begin{align}\label{vecsol}
\mathbf{s}_n(t)
=
e^{\mathbf{M}_n t} \mathbf{s}_n(0)
-
\mathbf{M}_n^{-1} \left(\mathbf{I} - e^{\mathbf{M}_n t}\right) \mathbf{c}_n ,
\end{align}
and if all $m_{k,k} \ne 0$ for all $k \in \{1, \dots, n\}$ are distinct, the $n^\text{th}$ scalar function $s_n(t)$ is given by
\begin{align}\label{scalsol}
s_n(t)
&=
\mathbf{m}_{n}
\left(
\mathbf{M}_{n-1} - m_{n,n}\mathbf{I}
\right)^{-1}
\left(
e^{\mathbf{M}_{n-1} t} - e^{m_{n,n}t}\mathbf{I}
\right)
\left(
\mathbf{s}_{n-1}(0)
+
\frac{\mathbf{c}_{n-1}}{m_{n,n}}
\right)
+
e^{m_{n,n}t}s_n(0)
\nonumber
\\
&
\quad
-
\frac{c_n}{m_{n,n}} \left( 1-e^{m_{n,n}t} \right)
+
\frac{\mathbf{m}_{n}}{m_{n,n}}\mathbf{M}_{n-1}^{-1}
\left(
\mathbf{I}
-
e^{\mathbf{M}_{n-1}t}
\right)
\mathbf{c}_{n-1}
,
\end{align}
where $t \geq 0$.
\end{lemma}

With this lemma in hand, we can now move to using these matrix sequences for calculating Markov process moments. To do so, we will use the  infinitesimal generator, a key tool for Markov processes, to find the derivatives of the moments through time. By identifying a Matryoshkan matrix structure in these differential equations, we are able to apply Lemma~\ref{mmlemma} to find closed form expressions for the moments.

\section{Calculating Moments through Matryoshkan Matrix Sequences}\label{secMoments}

In this section we connect Matryoshkan matrix sequences with the moments of Markov processes. \edit{In Subsection~\ref{genSec}, we prove the main result, which is the computation of the moment of a general Markov process through Matryoshkan matrices.} To demonstrate the applicability of this result, we now apply it to a series of examples. First in Subsection~\ref{subsecHawkes} we obtain the moments of the self-exciting Hawkes process, for which finding moments in closed form has been an open problem. Then in Subsection~\ref{subsecSN} we study the Markovian shot noise process, a stochastic intensity process that trades self-excitement for external shocks. Next in Subsection~\ref{subsecIto} we showcase the use of these techniques for diffusive dynamics through application to It\^o diffusions. Finally, in Subsection~\ref{subsecAQH} we apply this technique to a process with jumps both upwards and downwards, a linear birth-death-immigration process we have studied previously called the Affine Queue-Hawkes process. In each scenario, we describe the process of interest, define the infinitesimal generator, and identify the matrix structure. Through this, we solve for the process moments.

\subsection{The Moments of General Markov Processes}\label{genSec}

The connection between Matryoshkan matrices and Markov processes is built upon a key tool for Markov processes, the infinitesimal generator. \edit{For a detailed introduction to infinitesimal generators and their use in studying Markov processes, see e.g.~\citet{ethier2009markov}.}  For a Markov process $X_t$ on state space $\mathbb{S}$, the infinitesimal generator on a function $f : \mathbb{S} \to \mathbb{R}$ is defined
$$
\mathcal{L}f(x)
=
\lim_{\tau \to 0}
\frac{\E{f(X_\tau) \mid X_0 = x} - f(x)}{\tau}
.
$$
In our context and in many others, the power of the infinitesimal generator comes through use of Dynkin's formula, which gives us that
$$
\frac{\mathrm{d}}{\mathrm{d}t}\E{f(X_t)}
=
\E{\mathcal{L}f(X_t)}
.
$$
To study the moments of a Markov process, we are interested in functions $f$ that are polynomials.  Let's suppose now that $\mathcal{L}x^n$ for any $n \in \mathbb{Z}^+$ is polynomial in the lower powers of $x$ for a given Markov process $X_t$. Then, we can then write
$$
\mathcal{L}X_t^n
=
\theta_{0,n}
+
\sum_{i=1}^n \theta_{i,n} X_t^i
,
$$
which implies that the differential equation for the $n^\text{th}$ moment of this process is
$$
\frac{\mathrm{d}}{\mathrm{d}t} \E{X_t^n}
=
\theta_{0,n} + \sum_{i=1}^n \theta_{i,n} \E{X_t^i}
,
$$
for some collection of constants $\theta_{0,n}$, $\theta_{1,n}$, \dots, $\theta_{n,n}$. Thus, to solve for the $n^\text{th}$ moment of $X_t$ we must first solve for all the moments of lower order. We can also observe that to solve for the $(n-1)^\text{th}$ moment we must have all the moments below it. In comparing these systems of differential equations, we can see that all of the equations in the system for the $(n-1)^\text{th}$ moment are also in the system for the $n^\text{th}$ moment. No coefficients are changed in any of these lower moment equations, the only difference between the two systems is the inclusion of the differential equation for the $n^\text{th}$ moment in its own system. Hence, the nesting Matryoshkan structure arises. By expressing each system of linear ordinary differential equations in terms of a vector of moments, a matrix of coefficients, and a vector of shift terms, we can use these matrix sequences to capture how one moment's system encapsulates all the systems below it. This observation then allows us to calculate all the moments of the process in closed form, as we now show in Theorem~\ref{singlesol}.

\begin{theorem}\label{singlesol}
Let $X_t$ be a Markov process such that the time derivative of its $n^\text{th}$ moment can be written
\begin{align}
\frac{\mathrm{d}}{\mathrm{d}t} \E{X_t^n}
&=
\theta_{0,n} + \sum_{i=1}^n \theta_{i,n} \E{X_t^i}
,
\label{momentODE}
\end{align}
for any $n \in \mathbb{Z}^+$, where $t \geq 0$ and $\theta_{i, n} \in \mathbb{R}$ for all $i \leq n$. Let $\boldsymbol{\Theta}_n \in \mathbb{R}^{n \times n}$ be defined recursively by
\begin{align}
\mathbf{\Theta}_n
=
\begin{bmatrix}
\boldsymbol{\Theta}_{n-1} & \mathbf{0}_{n-1 \times 1} \\
\boldsymbol{\theta}_n & \theta_{n,n}
\end{bmatrix}
,
\label{Thetadef}
\end{align}
where $\boldsymbol{\theta}_n = [\theta_{1, n}, \, \dots, \, \theta_{n-1, n}]$ and $\boldsymbol{\Theta}_1 = \theta_{1,1}$. Furthermore, let $\boldsymbol{\theta}_{0,n} = [\theta_{0,1}, \, \dots, \, \theta_{0,n}]^\T$. Then, if $\theta_{k,k} \ne 0$ for all $k \in \{1, \dots , n\}$ are distinct, the $n^\text{th}$ moment of $X_t$ can be expressed
\begin{align}\label{scalsol}
\E{X_t^n}
&=
\boldsymbol{\theta}_{n}
\left(
\boldsymbol{\Theta}_{n-1} - \theta_{n,n}\mathbf{I}
\right)^{-1}
\left(
e^{\boldsymbol{\Theta}_{n-1} t} - e^{\theta_{n,n}t}\mathbf{I}
\right)
\left(
\mathbf{x}_{n-1}(x_0)
+
\frac{\boldsymbol{\theta}_{0,n-1}}{\theta_{n,n}}
\right)
+
x_0^n e^{\theta_{n,n}t}
-
\frac{\theta_{0,n}(1-e^{\theta_{n,n}t})}{\theta_{n,n}}
\nonumber
\\
&
\quad
+
\frac{\boldsymbol{\theta}_{n}}{\theta_{n,n}}\boldsymbol{\Theta}_{n-1}^{-1}
\left(
\mathbf{I}
-
e^{\boldsymbol{\Theta}_{n-1}t}
\right)
\boldsymbol{\theta}_{0,n-1}
,
\end{align}
where $x_0$ is the initial value of $X_t$ and where $\mathbf{x}_n(a) \in \mathbb{R}^n$ is such that $\left(\mathbf{x}_n(a)\right)_i = a^i$. If $X_t$ has a stationary distribution, then the $n^\text{th}$ steady-state moment $\E{X_\infty^n}$ is given by
\begin{align}\label{steadystate}
\E{X_\infty^n}
&=
\frac{1}{\theta_{n,n}}\left(\boldsymbol{\theta}_{n}\boldsymbol{\Theta}_{n-1}^{-1}\boldsymbol{\theta}_{0,n-1}
-
\theta_{0,n}\right)
,
\end{align}
and these steady-state moments satisfy the recursive relationship
\begin{align}
\E{X_\infty^{n+1}}
&=
-
\frac{1}{\theta_{n+1,n+1}}
\left(
\boldsymbol{\theta}_{n+1}
\mathbf{s}_{n}^X(\infty)
+
\theta_{0,n+1}
\right)
,
\end{align}
where $\mathbf{s}_n^X(\infty) \in \mathbb{R}^n$ is the vector of steady-state moments such that $\left(\mathbf{s}_n^X(\infty)\right)_i = \E{X_\infty^i}$.
\end{theorem}
\proof{Proof.}
Using the definition of $\boldsymbol{\Theta}_n$ in Equation~\ref{Thetadef}, Equation~\ref{momentODE} gives rise to the system of ordinary differential equations given by
\begin{align}\label{genDE}
\frac{\mathrm{d}}{\mathrm{d}t}\mathbf{s}_n^X(t)
=
\boldsymbol{\Theta}_n \mathbf{s}_n^X(t)
+
\boldsymbol{\theta}_{0,n}
,
\end{align}
where $\mathbf{s}_n^X(t) \in \mathbb{R}^n$ is the vector of transient moments at time $t \geq 0$ such that $\left(\mathbf{s}_n^X(t)\right)_i = \E{X_t^i}$ for all $1 \leq i \leq n$. We can observe that by definition the matrices $\boldsymbol{\Theta}_n$ form a Matryoshkan sequence, and thus by Lemma~\ref{mmlemma}, we achieve the stated transient solution. To prove the steady-state solution, we can first note that if the process has a steady-state distribution then the vector $\mathbf{s}_n^X(\infty) \in \mathbb{R}^n$ defined $\left(\mathbf{s}_n^X(\infty)\right)_i = \E{X_\infty^i}$ will satisfy
\begin{align}\label{sseq}
0
=
\boldsymbol{\Theta}_n \mathbf{s}_n^X(\infty)
+
\boldsymbol{\theta}_{0,n}
,
\end{align}
as this is the equilibrium solution to the differential equation corresponding to each of the moments. This system has a unique solution since $\boldsymbol{\Theta}_n$ is nonsingular due to the assumption that the diagonal values are unique and non-zero. Using Proposition~\ref{matryoshkanProp}, we find the $n^\text{th}$ moment by
$$
\E{X_\infty^n}
=
-\mathbf{v}_n^\T \boldsymbol{\Theta}_n^{-1} \boldsymbol{\theta}_{0,n}
=
\begin{bmatrix}
\frac{1}{\theta_{n,n}}\boldsymbol{\theta}_{n}\boldsymbol{\Theta}_{n-1}^{-1} & -\frac{1}{\theta_{n,n}}
\end{bmatrix}
\begin{bmatrix}
\boldsymbol{\theta}_{0,n-1} \\
\theta_{0,n}
\end{bmatrix}
=
\frac{1}{\theta_{n,n}}\left(\boldsymbol{\theta}_{n}\boldsymbol{\Theta}_{n-1}^{-1}\boldsymbol{\theta}_{0,n-1}
-
\theta_{0,n}\right)
,
$$
which completes the proof of Equation~\ref{steadystate}. To conclude, one can note that each line of the linear system in Equation~\ref{sseq} implies the stated recursion.
\Halmos\\
\endproof

%
%

\subsection{Application to Hawkes Process Intensities}\label{subsecHawkes}

For our first example of this method let us turn to our motivating application, the \edit{Markovian} Hawkes process intensity. Via \citet{hawkes1971spectra}, this process is defined as follows. Let $\lambda_t$ be stochastic arrival process intensity such that
\edit{
$$
\lambda_t
=
\lambda^* + (\lambda_0 - \lambda^*) e^{-\beta t} + \int_0^t \alpha e^{-\beta (t-s)} \mathrm{d}N_s
=
\lambda^* + (\lambda_0 - \lambda^*) e^{-\beta t} + \sum_{i=1}^{N_t} \alpha e^{-\beta (t - A_i)}
,
$$
}where $\{A_i \mid i \in \mathbb{Z}^+\}$ is the sequence of arrival epochs \edit{in the point process $N_t$}, with
$$
\PP{N_{t + s} - N_t = 0 \mid \mathcal{F}_t}
=
\PP{N_{t + s} - N_t = 0 \mid \lambda_t}
=
e^{-\int_0^s \lambda_{t + u} \mathrm{d}u}
,
$$
where $\mathcal{F}_t$ is the filtration generated by the history of $\lambda_t$ up to time $t$. We will assume that $\beta > \alpha > 0$ so that the process has a stationary distribution, and we will also let $\lambda^* > 0$ \edit{and $\lambda_0 > 0$}. Note that the process behaves as follows: at arrivals $\lambda_t$ increases by $\alpha$ and in the interims it decays exponentially at rate $\beta$ towards the baseline level $\lambda^*$. In this way, $(\lambda_t, N_t)$ is referred to as a \textit{self-exciting} point process, as the occurrence of an arrival increases the intensity and thus increases the likelihood that another arrival will occur soon afterwards. Because the intensity $\lambda_t$ forms a Markov process, we can write its infinitesimal generator for a (sufficiently regular) function $f:\mathbb{R}^+ \to \mathbb{R}$ as follows:
$$
\mathcal{L}f(\lambda_t) = \lambda_t \left(f(\lambda_t + \alpha) - f(\lambda_t)\right) - \beta  \left(\lambda_t -  \lambda^*\right) \frac{\mathrm{d}f(\lambda_t)}{\mathrm{d}\lambda_t}
.
$$
Note that this expression showcases the process dynamics that we have described, as the first term on the right-hand side corresponds to the product of the arrival rate and the change in the process when an arrival occurs while the second term captures the decay.

To obtain the $n^\text{th}$ moment we must consider $f(\cdot)$ of the form $f(x) = x^n$. In the simplest case when $n = 1$ this formula yields an ordinary differential equation for the mean, which can be written
$$
\frac{\mathrm{d}}{\mathrm{d}t}\E{\lambda_t}
=
\alpha\E{\lambda_t} - \beta\left(\E{\lambda_t} - \lambda^*\right)
=
\beta \lambda^*
-
(\beta - \alpha) \E{\lambda_t}
.
$$
By comparison for the second moment at $n = 2$ we have
$$
\frac{\mathrm{d}}{\mathrm{d}t}\E{\lambda_t^2}
=
\E{\lambda_t \left((\lambda_t + \alpha)^2 - \lambda_t^2\right) - 2 \beta \lambda_t (\lambda_t - \lambda^*)}
=
(2 \beta \lambda^* + \alpha^2)\E{ \lambda_t  }
-
2(\beta-\alpha)\E{\lambda_t^2}
,
$$
and we can note that while the ODE for the mean is autonomous, the second moment equation depends on both the mean and the second moment. Thus, to solve for the second moment we must also solve for the mean, leading us to the following system of differential equations:
$$
\frac{\mathrm{d}}{\mathrm{d}t}
\begin{bmatrix}
\E{\lambda_t}\\ \E{\lambda_t^2}
\end{bmatrix}
=
\begin{bmatrix}
-(\beta - \alpha) & 0 \\
2\beta\lambda^* + \alpha^2 & -2(\beta-\alpha)
\end{bmatrix}
\begin{bmatrix}
\E{\lambda_t}\\ \E{\lambda_t^2}
\end{bmatrix}
+
\begin{bmatrix}
\beta\lambda^* \\ 0
\end{bmatrix}
.
$$
Moving on to the third moment, the infinitesimal generator formula yields
$$
\frac{\mathrm{d}}{\mathrm{d}t}\E{\lambda_t^3}
=
\E{\lambda_t\left((\lambda_t + \alpha)^3 - \lambda_t^3\right) - 3\beta \lambda_t^{2}\left(\lambda_t - \lambda^*\right)}
=
\alpha^{3}\E{\lambda_t }
+
3(\beta\lambda^* + \alpha^2)\E{\lambda_t^{2}}
-
3(\beta - \alpha) \E{\lambda_t^{3} }
,
$$
and hence we see that this ODE now depends on all of the first three moments. Thus, to solve for $\E{\lambda_t^3}$ we need to solve the system of ordinary differential equations
$$
\frac{\mathrm{d}}{\mathrm{d}t}
\begin{bmatrix}
\E{\lambda_t}\\ \E{\lambda_t^2}\\ \E{\lambda_t^3}
\end{bmatrix}
=
\begin{bmatrix}
-(\beta - \alpha) & 0  & 0\\
2\beta\lambda^* + \alpha^2 & -2(\beta-\alpha) & 0 \\
\alpha^3 & 3(\beta\lambda^* + \alpha^2) & - 3(\beta-\alpha)
\end{bmatrix}
\begin{bmatrix}
\E{\lambda_t}\\ \E{\lambda_t^2}\\ \E{\lambda_t^3}
\end{bmatrix}
+
\begin{bmatrix}
\beta\lambda^* \\ 0 \\ 0
\end{bmatrix}
,
$$
and this now suggests the Matryoshkan structure of these process moments: we can note that the system for the second moment is nested within the system for the third moment. That is, the matrix for the three dimensional system contains the two dimensional system in its upper left-hand block, just as the vector of the first three moments has the first two moments in its first two coordinates. In general, we can see that the $n^\text{th}$ moment will satisfy the ODE given by
\begin{align*}
\frac{\mathrm{d}}{\mathrm{d}t}\E{\lambda_t^n}
&=
\E{\lambda_t\left((\lambda_t + \alpha)^n - \lambda_t^n\right) - n\beta \lambda_t^{n-1}\left(\lambda_t - \lambda^*\right)}
\\&=
\sum_{k=1}^{n}{n \choose k - 1}\alpha^{n-k+1}\E{\lambda_t^{k} }
-
n\beta \E{\lambda_t^{n} }
+
n\beta\lambda^* \E{\lambda_t^{n-1}}
,
\end{align*}
where we have simplified by use of the binomial theorem. Thus, the system of differential equations needed to solve for the $n^\text{th}$ moment uses the matrix from the $(n-1)^\text{th}$ system augmented below by the row
\begin{align}\label{augHawkes}
\begin{bmatrix}
\alpha^n & n\alpha^{n-1} & {n \choose 2}\alpha^{n-2} & \dots & {n \choose n - 3} \alpha^3 & {n \choose n - 2}\alpha^2 + n \beta\lambda^* & -n(\beta-\alpha)
\end{bmatrix}
,
\end{align}
and buffered on the right by a column of zeros. To collect these coefficients into a coherent structure, let us define the matrix $\boldsymbol{\mathcal{P}}_n(a) \in \mathbb{R}^{n \times n}$ for $a \in \mathbb{R}$ such that
\begin{align}\label{MPMdef}
\left(\boldsymbol{\mathcal{P}}_n(a)\right)_{i,j}
=
\begin{cases}
{i \choose j-1} a^{i-j+1} & i \geq j, \\
0 & i < j .
\end{cases}
\end{align}
If we momentarily disregard the terms with $\beta$ in the general augment row in Equation~\ref{augHawkes}, one can observe that the remaining terms in this vector are given by the bottom row of the matrix $\boldsymbol{\mathcal{P}}_n(\alpha)$. Furthermore, by definition $\{\boldsymbol{\mathcal{P}}_n(a) \mid n \in \mathbb{Z}^+\}$ forms a Matryoshkan matrix sequence. We can also note that $\boldsymbol{\mathcal{P}}_n(a)$ can be equivalently defined as
$$
\boldsymbol{\mathcal{P}}_n(a)
=
\sum_{k=1}^n
a
\begin{bmatrix}
\mathbf{0}_{n-k \times n-k} & \mathbf{0}_{n-k \times k}\\
\mathbf{0}_{k \times n-k} & \mathbf{L}_k(a)
\end{bmatrix}
,
$$
where $\mathbf{L}_k(a) = e^{a\mathbf{diag}\left(1:k-1,-1\right)}$ is the $k^\text{th}$ lower triangular Pascal matrix, i.e. the nonzero terms in $\mathbf{L}_k(1)$ yield the first $k$ rows of Pascal's triangle. Alternatively, $\boldsymbol{\mathcal{P}}_n(a)$ can be found by creating a lower triangular matrix from the strictly lower triangular values in $\mathbf{L}_{n+1}(a)$. For these reasons, we refer to the sequence of $\boldsymbol{\mathcal{P}}_n(a)$ as \textit{Matryoshkan Pascal matrices}. For brief overviews and beautiful properties of Pascal matrices, see \citet{brawer1992linear,call1993pascal,zhang1997linear,edelman2004pascal}. As we have seen in the preceding derivation, Matryoshkan Pascal matrices arise naturally in using the infinitesimal generator for calculating moments of Markov processes. This follows from the application of the binomial theorem to jump terms. Now, in the case of the Markovian Hawkes process intensity we find closed form expressions for all transient moments in Corollary~\ref{hawkesCor}.

\begin{corollary}\label{hawkesCor}
Let $\lambda_t$ be the intensity of a Hawkes process with baseline intensity $\lambda^* > 0$, intensity jump $\alpha > 0$, and decay rate $\beta > \alpha$. Then, the $n^\text{th}$ moment of $\lambda_t$ is given by
\begin{align}
\E{\lambda_t^n}
&=
\mathbf{m}^\lambda_{n}
\left(
\mathbf{M}^\lambda_{n-1} + n(\beta-\alpha)\mathbf{I}
\right)^{-1}
\left(
e^{\mathbf{M}^\lambda_{n-1} t} - e^{-n(\beta-\alpha)t}\mathbf{I}
\right)
\left(
\mathbf{x}_{n-1}(\lambda_0)
-
\frac{\beta \lambda^* \mathbf{v}_1 }{n(\beta-\alpha)}
\right)
+
\lambda_0^n  e^{-n(\beta-\alpha)t}
\nonumber
\\
&
\quad
+
\mathbb{I}_{\{n = 1\}}\frac{\beta\lambda^*}{\beta - \alpha}
\left(
1 - e^{-(\beta - \alpha) t}
\right)
-
\frac{\beta\lambda^*}{n(\beta-\alpha)}\mathbf{m}^\lambda_{n}\left(\mathbf{M}^\lambda_{n-1}\right)^{-1}
\left(
\mathbf{I}
-
e^{\mathbf{M}^\lambda_{n-1} t}
\right)
\mathbf{v}_{1}
,
\end{align}
for all $t \geq 0$ and $n \in \mathbb{Z}^+$, where \edit{$\mathbf{v}_1 \in \mathbb{R}^n$ is the unit vector in the first coordinate,} $\mathbf{M}_n^\lambda = \beta \lambda^* \mathbf{diag}\left(2:n,-1\right) - \beta \mathbf{diag}\left(1:n\right) + \boldsymbol{\mathcal{P}}_n(\alpha)$, $
\mathbf{m}_n^{\lambda}
=
\left[
\left(\mathbf{M}_n^\lambda\right)_{n,1}
,
\,
\dots
,
\,
\left(\mathbf{M}_n^\lambda\right)_{n,n-1}
\right]
$
is given by
\begin{align*}
\left(\mathbf{m}_n^{\lambda}\right)_{j}
&
=
\begin{cases}
{n \choose j-1} \alpha^{n-j+1}
&
\text{ if $j < n-1$,}
\\
{n \choose n-2} \alpha^{2} + n\beta \lambda^*
&
\text{ if $ j =n-1$,}
\end{cases}
\end{align*}
and $\mathbf{x}_n(a) \in \mathbb{R}^n$ is such that $\left(\mathbf{x}_n(a)\right)_i = a^i$. In steady-state, the $n^\text{th}$ moment of $\lambda_t$ is given by
\begin{align}
\lim_{t \to \infty}\E{\lambda_t^n}
=
-
\frac{\beta\lambda}{n(\beta - \alpha)}
\mathbf{m}_n^\lambda \left(\mathbf{M}_{n-1}^\lambda\right)^{-1}\mathbf{v}_1
,
\end{align}
for $n \geq 2$ with $\lim_{t \to \infty} \E{\lambda_t} = \frac{\beta \lambda^*}{\beta - \alpha}$. Moreover, the $(n+1)^\text{th}$ steady-state moment of the Hawkes process intensity is given by the recursion
\begin{align}\label{hawkesSSR}
\lim_{t \to \infty}\E{\lambda_t^{n+1}} = \frac{1}{(n+1)(\beta-\alpha)}\mathbf{m}_{n+1}^\lambda \mathbf{s}_n^\lambda ,
\end{align}
for all $n \in \mathbb{Z}^+$, where $\mathbf{s}_n^\lambda \in \mathbb{R}^n$ is the vector of steady-state moments defined such that $\left(\mathbf{s}_n^\lambda\right)_i = \lim_{t \to \infty}\E{\lambda_t^i}$ for $1 \leq i \leq n$.
\end{corollary}


\subsection{Application to Shot Noise Processes}\label{subsecSN}

As a second example of calculating moments through Matryoshkan matrices, consider a \edit{Markovian} shot noise process; \edit{see e.g.~\citet{daley2003models} for an introduction}. That is, let $\psi_t$ be defined such that
$$
\psi_t
=
\sum_{i=1}^{N_t} J_i e^{-\beta (t - A_i)}
,
$$
where $\beta > 0$, $\{J_i \mid i \in \mathbb{Z}^+\}$ is a sequence of i.i.d. positive random variables, $N_t$ is a Poisson process at rate $\lambda > 0$, and $\{A_i \mid i \in \mathbb{Z}^+\}$ is the sequence of arrival times in the Poisson process. These dynamics yield the following infinitesimal generator:
$$
\mathcal{L}f(\psi_t) = \lambda \left( f(\psi_t + J_i) - f(\psi_t)\right) - \beta \psi_t \frac{\mathrm{d}f(\psi_t)}{\mathrm{d}\psi_t}
.
$$
We can note that this is similar to the Hawkes process discussed in Subsection~\ref{subsecHawkes}, as the right-hand side contains a term for jumps and a term for exponential decay. However, this infinitesimal generator formula also shows key differences between the two processes, as the jumps in the shot noise process are of random size and they occur at the fixed, exogenous rate $\lambda > 0$. Supposing the mean jump size is finite, this now yields that the mean satisfies the ordinary differential equation
$$
\frac{\mathrm{d}}{\mathrm{d}t}\E{\psi_t}
=
\lambda \E{J_1} - \beta \E{\psi_t}
,
$$
whereas if $\E{J_1^2} < \infty$, the second moment of the shot noise process is given by the solution to
$$
\frac{\mathrm{d}}{\mathrm{d}t}\E{\psi_t^2}
=
\E{\lambda \left( (\psi_t + J_1)^2 - \psi_t^2 \right) - 2 \beta \psi_t^2}
=
\lambda \E{J_1^2} + 2\lambda \E{J_1} \E{\psi_t} - 2 \beta \E{\psi_t^2}
,
$$
which depends on both the second moment and the mean. This gives rise to the linear system of differential equations
$$
\frac{\mathrm{d}}{\mathrm{d}t}
\begin{bmatrix}
\E{\psi_t}\\ \E{\psi_t^2}
\end{bmatrix}
=
\begin{bmatrix}
-\beta & 0 \\
2\lambda\E{J_1} & -2\beta
\end{bmatrix}
\begin{bmatrix}
\E{\psi_t}\\ \E{\psi_t^2}
\end{bmatrix}
+
\begin{bmatrix}
\lambda \E{J_1} \\ \lambda \E{J_1^2}
\end{bmatrix}
,
$$
and by observing that the differential equation for the third moment depends on the first three moments if the third moment of the jump size is finite,
$$
\frac{\mathrm{d}}{\mathrm{d}t}\E{\psi_t^3}
=
\E{\lambda \left( (\psi_t + J_1)^3 - \psi_t^3 \right) - 3 \beta \psi_t^3}
=
\lambda \E{J_1^3} + 3\lambda \E{J_1^2} \E{\psi_t} + 3\lambda\E{J_1} \E{\psi_t^2} - 3 \beta \E{\psi_t^3}
,
$$
we can see that the system for the first two moments again are contained in the system for the first three moments:
$$
\frac{\mathrm{d}}{\mathrm{d}t}
\begin{bmatrix}
\E{\psi_t}\\ \E{\psi_t^2}\\ \E{\psi_t^3}
\end{bmatrix}
=
\begin{bmatrix}
-\beta & 0 & 0 \\
2\lambda\E{J_1} & -2\beta & 0 \\
3\lambda\E{J_1^2} & 3\lambda\E{J_1} & - 3\beta
\end{bmatrix}
\begin{bmatrix}
\E{\psi_t}\\ \E{\psi_t^2}\\ \E{\psi_t^3}
\end{bmatrix}
+
\begin{bmatrix}
\lambda \E{J_1} \\ \lambda \E{J_1^2} \\ \lambda \E{J_1^3}
\end{bmatrix}
.
$$
By use of the binomial theorem, we can observe that if $\E{J_1^n} < \infty$ then the $n^\text{th}$ moment of the shot noise process satisfies
$$
\frac{\mathrm{d}}{\mathrm{d}t}\E{\psi_t^n}
=
\E{\lambda \left( (\psi_t + J_1)^n - \psi_t^n \right) - n \beta \psi_t^n}
=
\sum_{k=0}^{n-1} {n \choose k}\E{J_1^{n-k}}\E{\psi_t^k} - n \beta \E{\psi_t^n}
,
$$
which means that the $n^\text{th}$ dimensional system is equal to the preceding one augmented below by the row vector
$$
\begin{bmatrix}
n\lambda\E{J_1^{n-1}} & {n \choose 2}\lambda\E{J_1^{n-2}} & {n \choose 3}\lambda\E{J_1^{n-3}} & \dots & {n \choose n - 2}\lambda \E{J_1^2} & n\lambda\E{J_1} & -n\beta
\end{bmatrix}
,
$$
and to the right by zeros. Bringing this together, this now leads us to Corollary~\ref{SNCor}.

\begin{corollary}\label{SNCor}
Let $\psi_t$ be the intensity of a shot noise process with epochs given by a Poisson process with rate $\lambda > 0$, jump sizes drawn from the i.i.d. sequence of random variables $\{J_i \mid i \in \mathbb{Z}^+\}$, and exponential decay at rate $\beta > 0$. If $\E{J_1^n} < \infty$, the $n^\text{th}$ moment of $\psi_t$ is given by
\begin{align}
\E{\psi_t^n}
&=
\mathbf{m}^{\psi}_{n}
\left(
\mathbf{M}^{\psi}_{n-1} + n \beta\mathbf{I}
\right)^{-1}
\left(
e^{\mathbf{M}^{\psi}_{n-1} t} - e^{-n \beta t}\mathbf{I}
\right)
\left(
\mathbf{x}_{n-1}\left(\psi_0\right)
-
\frac{\lambda \mathbf{j}_{n-1}}{n\beta}
\right)
+
\psi_0^n e^{-n\beta t}
\nonumber
\\
&
\quad
+
\frac{\lambda\E{J_1^n}}{n\beta}\left(1-e^{-n\beta t}\right)
-
\frac{\lambda \mathbf{m}^{\psi}_{n}}{n\beta}\left(\mathbf{M}^{\psi}_{n-1}\right)^{-1}
\left(
\mathbf{I}
-
e^{\mathbf{M}^{\psi}_{n-1}t}
\right)
\mathbf{j}_{n-1}
,
\end{align}
for all $t \geq 0$ and $n \in \mathbb{Z}^+$, where $\mathbf{j}_{n} \in \mathbb{R}^n$ is such that $\left(\mathbf{j}_{n}\right)_i = \E{J_1^i}$, $\mathbf{M}_n^\psi \in \mathbb{R}^{n \times n}$ is recursively defined
$$
\mathbf{M}_n^\psi
=
\begin{bmatrix}
\mathbf{M}_{n-1}^\psi & \mathbf{0}_{n-1 \times 1} \\
\mathbf{m}_n^\psi & - n\beta
\end{bmatrix}
,
$$
with the row vector $\mathbf{m}_n^{\psi} \in \mathbb{R}^{n-1}$ defined such that $\left(\mathbf{m}_n^{\psi}\right)_i = {n \choose i}\lambda \E{J_1^{n-i}}$ and with $\mathbf{M}_1^\psi = -\beta$,
and where $\mathbf{x}_n(a) \in \mathbb{R}^n$ is such that $\left(\mathbf{x}_n(a)\right)_i = a^i$. In steady-state, the $(n+1)^\text{th}$ moment of the shot noise process is given by
\begin{align}
\lim_{t \to \infty} \E{\psi_t^n}
&=
\frac{\lambda}{n\beta}\left(\E{J_1^n} - \mathbf{m}_n^\psi \left( \mathbf{M}_{n-1}^\psi\right)^{-1} \mathbf{j}_{n-1} \right)
,
\end{align}
for $n \geq 2$ where $\lim_{t \to \infty} \E{\psi_t} = \frac{\lambda}{\beta}\E{J_1}$. Moreover, if $\E{J_1^{n+1}} < \infty$ the $(n+1)^\text{th}$ moment of the shot noise process is given by the recursion
\begin{align}
\lim_{t \to \infty} \E{\psi_t^{n+1}}
&=
\frac{1}{(n+1)\beta}\left(\mathbf{m}_{n+1}^\psi \mathbf{s}_n^\psi + \E{J_1^{n+1}}\right)
,
\end{align}
for all $n \in \mathbb{Z}^+$, where $\mathbf{s}_n^\psi \in \mathbb{R}^n$ is the vector of steady-state moments defined such that $(\mathbf{s}_n^\psi)_i = \lim_{t \to \infty} \E{\psi_t^i}$ for $1 \leq i \leq n$.
\end{corollary}

\subsection{Application to It\^{o} Diffusions}\label{subsecIto}

For our third example, we consider an It\^{o} diffusion; \edit{see e.g.~\citet{oksendal2013stochastic} for an overview}. Let $S_t$ be given by the stochastic differential equation
$$
\mathrm{d}S_t = g(S_t) \mathrm{d}t + h(S_t) \mathrm{d}B_t
,
$$
where $B_t$ is a Brownian motion and $g(\cdot)$ and $h(\cdot)$ are real-valued functions. Then, the infinitesimal generator for this process is given by
$$
\mathcal{L}f(S_t)
=
g(S_t)\frac{\mathrm{d}f(S_t)}{\mathrm{d}S_t}
+
\frac{h(S_t)^2}{2} \frac{\mathrm{d}^2f(S_t)}{\mathrm{d}S_t^2}
,
$$
where $f : \mathbb{R} \to \mathbb{R}$. Because we will be considering functions of the form $f(x) = x^n$ for $n \in \mathbb{Z}^+$, we will now specify the forms of $g(\cdot)$ and $h(\cdot)$ to be $g(x) = \mu + \theta x$ for some $\mu \in \mathbb{R}$ and $\theta \in \mathbb{R}$ and $h(x) = \sigma x^{\gamma\slash 2}$ for some $\sigma \in \mathbb{R}$ and $\gamma \in \{0,1,2\}$. One can note that this encapsulates a myriad of relevant stochastic processes including many that are popular in the financial models literature, such as Ornstein-Uhlenbeck (OU) processes, geometric Brownian motion (GBM), and Cox-Ingersoll-Ross (CIR) processes. In this case, the infinitesimal generator becomes
$$
\mathcal{L}f(S_t)
=
\left(\mu + \theta S_t\right)\frac{\mathrm{d}f(S_t)}{\mathrm{d}S_t}
+
\frac{\sigma^2 S_t^{\gamma}}{2} \frac{\mathrm{d}^2f(S_t)}{\mathrm{d}S_t^2}
,
$$
meaning that we can express the ordinary differential equation for the mean as
$$
\frac{\mathrm{d}}{\mathrm{d}t}\E{S_t}
=
\mu + \theta \E{S_t}
,
$$
and similarly the second moment will be given by the solution to
$$
\frac{\mathrm{d}}{\mathrm{d}t}\E{S_t^2}
=
\E{2(\mu + \theta S_t) S_t + \sigma^2 S_t^{\gamma}}
=
2\mu \E{S_t} + 2\theta \E{S_t^2} + \sigma^2 \E{S_t^\gamma}
.
$$
For the sake of example, we now let $\gamma = 1$ as is the case in the CIR process. Then, the first two transient moments of $S_t$ will be given by the solution to the system
$$
\frac{\mathrm{d}}{\mathrm{d}t}
\begin{bmatrix}
\E{S_t} \\
\E{S_t^2}
\end{bmatrix}
=
\begin{bmatrix}
\theta & 0 \\
2 \mu + \sigma^2 & 2\theta
\end{bmatrix}
\begin{bmatrix}
\E{S_t} \\
\E{S_t^2}
\end{bmatrix}
+
\begin{bmatrix}
\mu \\
0
\end{bmatrix}
.
$$
By observing that the third moment differential equation is
$$
\frac{\mathrm{d}}{\mathrm{d}t}\E{S_t^3}
=
\E{3(\mu + \theta S_t) S_t^2  + 3\sigma^2 S_t^{\gamma + 1}}
=
3\mu \E{S_t^2} + 3\theta \E{S_t^3} + 3 \sigma^2 \E{S_t^{\gamma + 1}}
,
$$
we can note that the third moment system for $\gamma = 1$ is
$$
\frac{\mathrm{d}}{\mathrm{d}t}
\begin{bmatrix}
\E{S_t} \\
\E{S_t^2} \\
\E{S_t^3}
\end{bmatrix}
=
\begin{bmatrix}
\theta & 0 & 0 \\
2 \mu + \sigma^2 & 2\theta & 0 \\
0 & 3\mu + 3 \sigma^2 & 3\theta
\end{bmatrix}
\begin{bmatrix}
\E{S_t} \\
\E{S_t^2} \\
\E{S_t^3}
\end{bmatrix}
+
\begin{bmatrix}
\mu \\
0 \\
0
\end{bmatrix}
,
$$
and this showcases the Matryoshkan nesting structure, as the second moment system is contained within the third. Because the general $n^\text{th}$ moment for $n \geq 2$ has differential equation given by
$$
\frac{\mathrm{d}}{\mathrm{d}t}\E{S_t^n}
=
\E{n(\mu + \theta S_t) S_t^{n-1}  + \frac{n(n-1)\sigma^2}{2} S_t^{n + \gamma - 2}}
=
n\mu \E{S_t^{n-1}} + n\theta \E{S_t^n} + \frac{n(n-1)\sigma^2}{2} \E{S_t^{n + \gamma - 2}}
,
$$
we can see that the $(n-1)^\text{th}$ system can be augmented below by the row vector
$\gamma = 1$
$$
\begin{bmatrix}
0 & 0 & \dots & 0 & n\mu +  \frac{n(n-1)\sigma^2}{2} & n\theta
\end{bmatrix}
,
$$
and to the right by zeros. Through this observation, we can now give the moments of It\^o diffusions in Corollary~\ref{itoCor}.

\begin{corollary}\label{itoCor}
Let $S_t$ be an It\^{o} diffusion that satisfies the stochastic differential equation
\begin{align}
\mathrm{d}S_t
=
(\mu + \theta S_t)\mathrm{d}t
+
\sigma S_t^{\gamma \slash 2} \mathrm{d} B_t
,
\end{align}
where $B_t$ is a Brownian motion and with $\mu, \theta, \sigma \in \mathbb{R}$ and $\gamma \in \{0, 1, 2\}$. Then, the $n^\text{th}$ moment of $S_t$ is given by
\begin{align}
\E{S_t^n}
&=
\mathbf{m}^S_{n}
\left(
\mathbf{M}^S_{n-1} - \chi_n\mathbf{I}
\right)^{-1}
\left(
e^{\mathbf{M}^S_{n-1} t} - e^{\chi_n t}\mathbf{I}
\right)
\left(
\mathbf{x}_{n-1}(S_0)
+
\frac{\mu \mathbf{v}_1 + \sigma^2\mathbb{I}_{\{\gamma = 0\}}\mathbf{v}_2}{\chi_n}
\right)
+
S_0^n e^{\chi_n t}
\nonumber
\\
&
\quad
-
\left(\mu \mathbb{I}_{\{n=1\}} + \sigma^2\mathbb{I}_{\{\gamma = 0, n= 2\}}\right)
\frac{1-e^{\chi_n t}}{\chi_n}
+
\frac{\mathbf{m}^S_{n}}{\chi_{n}}\left(\mathbf{M}^S_{n-1} \right)^{-1}
\left(
\mathbf{I}
-
e^{\mathbf{M}^S_{n-1}t}
\right)
\left(\mu\mathbf{v}_1 + \sigma^2\mathbb{I}_{\{\gamma = 0\}}\mathbf{v}_2\right)
,
\end{align}
for all $t \geq 0$ and $n \in \mathbb{Z}^+$, where $\chi_n = n\theta + \frac{n}{2}(n-1)\sigma^2\mathbb{I}_{\{\gamma = 2\}}$,
$\mathbf{M}^S_n = \theta \mathbf{diag}(1:n) + \mu \mathbf{diag}(2:n,-1) + \frac{\sigma^2}{2}\mathbf{diag}(\mathbf{d}^{2-\gamma}_{n+\gamma-2}, \gamma - 2)$ for $\mathbf{d}^j_k \in \mathbb{R}^k$ such that $(\mathbf{d}_k^j)_i = (j+i)(j+i-1)$, and
$\mathbf{m}_n^{S}
=
\left[
\left(\mathbf{M}_n^S\right)_{n,1}
,
\,
\dots
,
\,
\left(\mathbf{M}_n^S\right)_{n,n-1}
\right]
$
is such that
$$
\left(\mathbf{m}_n^{S}\right)_{j}
=
\begin{cases}
n\mu + \frac{n(n-1)\sigma^2}{2}\mathbb{I}_{\{\gamma = 1\}} & j = n -1, \\
\frac{n(n-1)\sigma^2}{2}\mathbb{I}_{\{\gamma = 0\}} & j = n - 2, \\
0 & 1 \leq j < n - 2,
\end{cases}
$$
and  $\mathbf{x}_n(a) \in \mathbb{R}^n$ is such that $\left(\mathbf{x}_n(a)\right)_i = a^i$. If \edit{$\theta < 0$ and $\gamma \in \{0,1\}$}, then the $n^\text{th}$ steady-state moment of $S_t$ is given by
\begin{align}
\lim_{t \to \infty}
\E{S_t^n}
=
\frac{\mu}{\chi_n}
\mathbf{m}_n^S \left(\mathbf{M}_{n-1}^S\right)^{-1} \mathbf{v}_1
,
\end{align}
for $n \geq 2$ with $\lim_{t \to \infty} \E{S_t} = -\frac{\mu}{\theta}$. Moreover, the $(n+1)^\text{th}$ steady-state moment of $S_t$ is given by the recursion
\begin{align}
\lim_{t \to \infty} \E{S_t^{n+1}}
&=
-
\frac{1}{\chi_n}
\mathbf{m}_{n+1}^S \mathbf{s}_n^S
,
\end{align}
for all $n \in \mathbb{Z}^+$, where $\mathbf{s}_n^S \in \mathbb{R}^n$ is the vector of steady-state moments defined such that $\left(\mathbf{s}_n^S\right)_i = \lim_{t \to \infty}\E{S_t^i}$ for $1 \leq i \leq n$.
\end{corollary}

As a consequence of these expressions we can also gain insight for the moments of an It\^o diffusion in the case of non-integer $\gamma \in [0,2]$, as is used in volatility models such as the CEV model and the SABR model, see e.g. \citet{henry2008analysis}. This can be achieve through bounding the differential equations, as the $n^\text{th}$ moment of such a diffusion is again given by
$$
\frac{\mathrm{d}}{\mathrm{d}t}\E{S_t^n}
=
n\mu \E{S_t^{n-1}} + n\theta \E{S_t^n} + \frac{n(n-1)\sigma^2}{2} \E{S_t^{n + \gamma - 2}}
,
$$
and the right-most term in this expression can be bounded above and below
$$
 \E{S_t^{n + \lfloor\gamma\rfloor - 2}}
\leq
 \E{S_t^{n + \gamma - 2}}
\leq
 \E{S_t^{n + \lceil\gamma\rceil - 2}}
,
$$
and the differential equations given by substituting these bounded terms form a closed system solvable by Corollary~\ref{itoCor}. Assuming the true differential equation and the upper and lower bounds all share an initial value, the solution to the bounded equations bounds the solution to the true moment equation, see \citet{hale2013introduction}.

\subsection{Application to Growth-Collapse Processes}\label{gcpSubsec}

For a fourth example, we consider growth-collapse processes with Poisson driven shocks.  These processes have been studied in variety of contexts, see e.g. \citet{boxma2006markovian, kella2009growth, kella2010markov,boxma2011some}.  More recently, these processes and their related extensions have seen renewed interest in the study of the crypto-currency Bitcoin, see for example \citet{frolkova2019bitcoin,koops2018predicting,javier2019exact,fralix2019classes}.  While growth-collapse processes can be defined in many different ways, for this example we use a definition in the style of Section 4 from \citet{boxma2006markovian}. We let $Y_t$ be the state of the growth collapse model and let $\{U_i \mid i \in \mathbb{Z}^+\}$ be a sequence of independent \edit{$\mathrm{Uni}(0,1)$} random variables that are also independent from the state and history of the growth-collapse process. Then, the infinitesimal generator of $Y_t$ is given by
$$
\mathcal{L}f(Y_t)
=
 \lambda  \frac{\mathrm{d}f(Y_t)}{\mathrm{d}Y_t} +  \mu   \left( f( U_i  Y_t ) - f( Y_t ) \right)
.
$$
Thus, $Y_t$ experiences linear growth at rate $\lambda > 0$ throughout time but it also collapses at epochs given by a Poisson process with rate $\mu > 0$. At the $i^\text{th}$ collapse epoch the process falls to a fraction of its current level, specifically it jumps down to $U_i Y_t$. Using the infinitesimal generator, we can see that the mean of this growth-collapse process satisfies
$$
\frac{\mathrm{d}}{\mathrm{d}t}\E{Y_t}
=
\lambda + \mu \left( \E{U_1 Y_t} - \E{Y_t} \right)
=
\lambda - \frac{\mu}{2} \E{Y_t}
,
$$
and its second moment will satisfy
$$
\frac{\mathrm{d}}{\mathrm{d}t}\E{Y_t^2}
=
2 \lambda \E{Y_t}  + \mu \left( \E{U_1^2 Y_t^2} - \E{Y_t^2} \right)
=
2\lambda \E{Y_t} - \frac{2 \mu}{3} \E{Y_t^2}
.
$$
Therefore, we can write the linear system of differential equations for the second moment of this growth-collapse process as
$$
\frac{\mathrm{d}}{\mathrm{d}t}
\begin{bmatrix}
\E{Y_t}\\ \E{Y_t^2}
\end{bmatrix}
=
\begin{bmatrix}
-\frac{\mu}2 & 0 \\
2\lambda & - \frac{2\mu}{3}
\end{bmatrix}
\begin{bmatrix}
\E{Y_t}\\ \E{Y_t^2}
\end{bmatrix}
+
\begin{bmatrix}
\lambda \\  0
\end{bmatrix}
.
$$
Moving to the third moment, via the infinitesimal generator we write its differential equation as
$$
\frac{\mathrm{d}}{\mathrm{d}t}\E{Y_t^3}
=
3 \lambda \E{Y_t^2}  + \mu \left( \E{U_1^3 Y_t^3} - \E{Y_t^3} \right)
=
3\lambda \E{Y_t^2} - \frac{3 \mu}{4} \E{Y_t^3}
,
$$
which thus shows that the system of differential equations for the third moment is
$$
\frac{\mathrm{d}}{\mathrm{d}t}
\begin{bmatrix}
\E{Y_t}\\ \E{Y_t^2}\\ \E{Y_t^3}
\end{bmatrix}
=
\begin{bmatrix}
-\frac{\mu}2 & 0  & 0\\
2\lambda & - \frac{2\mu}{3} & 0\\
0 & 3\lambda & - \frac{3\mu}{4}
\end{bmatrix}
\begin{bmatrix}
\E{Y_t}\\ \E{Y_t^2}\\ \E{Y_t^3}
\end{bmatrix}
+
\begin{bmatrix}
\lambda \\  0\\ 0
\end{bmatrix}
,
$$
which evidently encapsulates the system for the first two moments. We can note that the general $n^\text{th}$ moment will satisfy
$$
\frac{\mathrm{d}}{\mathrm{d}t}\E{Y_t^n}
=
n \lambda \E{Y_t^{n-1}}  + \mu \left( \E{U_1^n Y_t^n} - \E{Y_t^n} \right)
=
n\lambda \E{Y_t^{n-1}} - \frac{n \mu}{n+1} \E{Y_t^n}
,
$$
and thus the system for the $n^\text{th}$ moment is given by appending the row vector
$$
\begin{bmatrix}
0 & 0 & \dots & 0 & n\lambda & -\frac{n\mu}{n+1}
\end{bmatrix}
$$
below the matrix from the $(n-1)^\text{th}$ system augmented by zeros on the right. Following this derivation, we reach the following general expressions for the moments in Corollary~\ref{gcpCor}. Furthermore, we can note that because of the relative simplicity of this particular structure, we are able to solve the recursion for the steady-state moments and give these terms explicitly.

\begin{corollary}\label{gcpCor}
Let $Y_t$ be a growth-collapse process with growth rate $\lambda > 0$ and uniformly sized collapses occurring according to a Poisson process with rate $\mu > 0$. Then, the $n^\text{th}$ moment of $Y_t$ is given by
\begin{align}
\E{Y_t^n}
&=
n\lambda\mathbf{v}_{n-1}^\T
\left(
\mathbf{M}^Y_{n-1} + \frac{n\mu}{n+1}\mathbf{I}
\right)^{-1}
\left(
e^{\mathbf{M}^Y_{n-1} t} - e^{-\frac{n\mu t}{n+1}}\mathbf{I}
\right)
\left(
\mathbf{x}_{n-1}(y_0)
-
\frac{(n+1)\lambda \mathbf{v}_1}{n\mu}
\right)
+
y_0^n e^{-\frac{n\mu t}{n+1}}
\nonumber
\\
&
\quad
+
\frac{(n+1)\lambda \mathbb{I}_{\{n=1\}}}{n\mu}
\left(1-e^{-\frac{n\mu t}{n+1}}\right)
-
\frac{(n+1)\lambda^2}{\mu}\mathbf{v}^\T_{n-1}\left(\mathbf{M}^Y_{n-1}\right)^{-1}
\left(
\mathbf{I}
-
e^{\mathbf{M}^Y_{n-1}t}
\right)
\mathbf{v}_1
,
\end{align}
where $y_0$ is the initial value of $Y_t$, $\mathbf{x}_n(a) \in \mathbb{R}^n$ is such that $\left(\mathbf{x}_n(a)\right)_i = a^i$, $\mathbf{M}_n^Y = \lambda \mathbf{diag}\left(2:n,-1\right) - \mu \mathbf{diag}\left(\frac{1}{2} : \frac{n}{n+1}\right)$, and $\mathbf{m}_n^Y = [\left(\mathbf{M}_n^Y\right)_{n,1}, \dots, \left(\mathbf{M}_n^Y\right)_{n,n-1}]$ is such that $\mathbf{m}_n^Y = n\lambda\mathbf{v}_{n-1}^\T$. Moreover the $n^\text{th}$ steady-state moment of $Y_t$ is given by
\begin{align}\label{steadystate}
\lim_{t \to \infty}\E{Y_t^n}
&=
2n!\left(\frac{\lambda}{\mu}\right)^n
,
\end{align}
for $n \in \mathbb{Z}^+$.
\end{corollary}

\subsection{Application to Ephemerally Self-Exciting Processes}\label{subsecAQH}

As a final detailed example of the applicability of Matryoshkan matrices, we now consider a stochastic process we have analyzed in \citet{daw2018queue}. This process is a linear birth-death-immigration process in which the occurrence of an arrival increases the arrival rate by an amount $\alpha > 0$, like in the Hawkes process, and this increase expires after an exponentially distributed duration with some rate $\mu > \alpha$. In this way, this process is an \textit{ephemerally self-exciting} process. Given a baseline intensity $\nu^* > 0$, let $Q_t$ be such that new arrivals occur at rate $\nu^* + \alpha Q_t$ and then the overall rate until the next excitement expiration is $\mu Q_t$. One can then think of $Q_t$ as the number of entities still causing active excitement at time $t \geq 0$. We will refer to $Q_t$ as the number in system for this ephemerally self-exciting process. The infinitesimal generator for a function $f : \mathbb{N}  \to \mathbb{R}$ is thus
$$
\mathcal{L}f(Q_t) = \left(\nu^* + \alpha Q_t\right) \left(f(Q_t + 1) - f(Q_t)\right) + \mu Q_t  \left(f(Q_t - 1) - f(Q_t)\right)
,
$$
which again captures the dynamics of the process, as the first term on the right hand-side is the product of the up-jump rate and the change in function value upon an increase in the process while the second term is the product of the down-jump rate and the corresponding process decrease. This now yields an ordinary differential equation for the mean given by
$$
\frac{\mathrm{d}}{\mathrm{d}t}\E{Q_t}
=
\nu^* + \alpha\E{Q_t} - \mu \E{Q_t}
=
\nu^*
-
(\mu - \alpha) \E{Q_t}
,
$$
while the second moment will satisfy
\begin{align*}
\frac{\mathrm{d}}{\mathrm{d}t}\E{Q_t^2}
&=
\E{\left(\nu^* + \alpha Q_t\right)\left((Q_t+1)^2 - Q_t^2\right) + \mu Q_t \left((Q_t-1)^2 - Q_t^2\right)}
\\
&=
\left(2\nu^* + \mu + \alpha\right)\E{Q_t}
+ \nu^*
- 2(\mu - \alpha)\E{Q_t^2}
.
\end{align*}
Thus, the first two moments are given by the solution to the linear system
$$
\frac{\mathrm{d}}{\mathrm{d}t}
\begin{bmatrix}
\E{Q_t}\\ \E{Q_t^2}
\end{bmatrix}
=
\begin{bmatrix}
-(\mu - \alpha) & 0 \\
2\nu^* + \mu + \alpha & -2(\mu-\alpha)
\end{bmatrix}
\begin{bmatrix}
\E{Q_t}\\ \E{Q_t^2}
\end{bmatrix}
+
\begin{bmatrix}
\nu^* \\  \nu^*
\end{bmatrix}
,
$$
and by observing that the third moment differential equation is
\begin{align*}
\frac{\mathrm{d}}{\mathrm{d}t}\E{Q_t^3}
&=
\E{\left(\nu^* + \alpha Q_t\right)\left((Q_t+1)^3 - Q_t^3\right) + \mu Q_t \left((Q_t-1)^3 - Q_t^3\right)}
\\
&=
\left(3\nu^* + 3\alpha + 3\mu\right)\E{Q_t^2}
+
\left(3\nu^* + \alpha - \mu\right)\E{Q_t}
+
\nu^*
- 3(\mu - \alpha)\E{Q_t^3}
,
\end{align*}
we can observe that the third moment system does indeed encapsulate that of the second moment:
$$
\frac{\mathrm{d}}{\mathrm{d}t}
\begin{bmatrix}
\E{Q_t}\\ \E{Q_t^2} \\ \E{Q_t^3}
\end{bmatrix}
=
\begin{bmatrix}
-(\mu - \alpha) & 0 & 0 \\
2\nu^* + \mu + \alpha & -2(\mu-\alpha) & 0 \\
3\nu^* + \alpha - \mu & 3\nu^* + 3\alpha + 3\mu & - 3(\mu-\alpha)
\end{bmatrix}
\begin{bmatrix}
\E{Q_t}\\ \E{Q_t^2} \\ \E{Q_t^3}
\end{bmatrix}
+
\begin{bmatrix}
\nu^* \\   \nu^* \\  \nu^*
\end{bmatrix}
.
$$
In general, the $n^\text{th}$ moment is given by the solution to
\begin{align*}
\frac{\mathrm{d}}{\mathrm{d}t}\E{Q_t^n}
&=
\E{(\nu^* + \alpha Q_t) \left((Q_t + 1)^n - Q_t^n\right) + \mu Q_t \left((Q_t - 1)^n - Q_t^n\right)}
\\
&=
\nu^*
+
\nu^* \sum_{k=1}^{n-1}{n \choose k} \E{Q_t^k}
+
\alpha \sum_{k=1}^{n}{n \choose k-1} \E{Q_t^k}
+
\mu \sum_{k=1}^{n} {n \choose k-1} \E{Q_t^{k}}(-1)^{n-k-1}
,
\end{align*}
which means that the $n^\text{th}$ system is given by augmenting the previous system below by
$$
\begin{bmatrix}
n\nu^* + \alpha + \mu(-1)^n
&
{n \choose 2}\nu^* + n\alpha + n\mu(-1)^{n-1}
&
\dots
&
n\nu^* + {n \choose 2}\alpha + {n \choose 2}\mu
&
-n(\mu-\alpha)
\end{bmatrix}
,
$$
and to the right by zeros. By comparing this row vector to the definition of the Martyoshkan Pascal matrices in Equation~\ref{MPMdef}, we arrive at explicit forms for the moments of this process shown now in Corollary~\ref{AQHCor}.

\begin{corollary}\label{AQHCor}
Let $Q_t$ be the number in system for an ephemerally self-exciting process with baseline intensity $\nu^* > 0$, intensity jump $\alpha > 0$, and duration rate $\mu > \alpha$. Then, the $n^\text{th}$ moment of $Q_t$ is given by
\begin{align}\label{scalsol}
\E{Q_t^n}
&=
\mathbf{m}^Q_{n}
\left(
\mathbf{M}^Q_{n-1} + n(\mu-\alpha)\mathbf{I}
\right)^{-1}
\left(
e^{\mathbf{M}^Q_{n-1} t} - e^{-n(\mu-\alpha)t}\mathbf{I}
\right)
\left(
\mathbf{x}_{n-1}(Q_0)
-
\frac{\nu^* \mathbf{v}}{n(\mu-\alpha)}
\right)
+
Q_0^n e^{-n(\mu-\alpha)t}
\nonumber
\\
&
\quad
+
\frac{\nu^*}{n(\mu-\alpha)}\left(1-e^{-n(\mu-\alpha)t}\right)
-
\frac{\nu^* \mathbf{m}^Q_{n}}{n(\mu-\alpha)}\left(\mathbf{M}^Q_{n-1}\right)^{-1}
\left(
\mathbf{I}
-
e^{\mathbf{M}^Q_{n-1}t}
\right)
\mathbf{v}
,
\end{align}
for all $t \geq 0$ and $n \in \mathbb{Z}^+$, where $\mathbf{M}^Q_n = \nu^*\boldsymbol{\mathcal{P}}_n(1)\mathbf{diag}(\mathbf{v},-1) + \alpha\boldsymbol{\mathcal{P}}_n(1) + \mu\boldsymbol{\mathcal{P}}_n(-1)$, and
$\mathbf{m}_n^{Q}
=
\left[
\left(\mathbf{M}_n^Q\right)_{n,1}
,
\,
\dots
,
\,
\left(\mathbf{M}_n^Q\right)_{n,n-1}
\right]
$
is such that
$
\left(\mathbf{m}_n^{Q}\right)_{j}
=
{n \choose j} \nu^*
+
{n \choose j-1}\alpha
+
{n \choose j-1}\mu (-1)^{n-j-1}
$
and  $\mathbf{x}_n(a) \in \mathbb{R}^n$ is such that $\left(\mathbf{x}_n(a)\right)_i = a^i$. In steady-state, the $n^\text{th}$ moment of $Q_t$ is given by
\begin{align}
\lim_{t \to \infty}
\E{Q_t^n}
=
\frac{\nu^*}{\mu - \alpha}
\left(
1 - \mathbf{m}_n^Q \left(\mathbf{M}_{n-1}^Q\right)^{-1} \mathbf{v}
\right)
,
\end{align}
for $n \geq 2$ with $\lim_{t \to \infty} \E{Q_t} = \frac{\nu^*}{\mu-\alpha}$. Moreover the $(n+1)^\text{th}$ steady-state moment of the ephemerally self-exciting process is given by the recursion
\begin{align}
\lim_{t \to \infty} \E{Q_t^{n+1}}
&=
\frac{1}{(n+1)(\mu-\alpha)}
\left(
\mathbf{m}_{n+1}^Q \mathbf{s}_n^Q + \nu^*
\right)
,
\end{align}
for all $n \in \mathbb{Z}^+$, where $\mathbf{s}_n^Q \in \mathbb{R}^n$ is the vector of steady-state moments defined such that $\left(\mathbf{s}_n^Q\right)_i = \lim_{t \to \infty}\E{Q_t^i}$ for $1 \leq i \leq n$.
\end{corollary}

\subsection{Additional Applications by Combination and Permutation}

While the preceding examples are the the only detailed examples we include in this paper, we can note that these Matryoshkan matrix methods can be applied to many other settings. In fact, one can observe that these example derivations can be applied directly to processes that feature a combination of their structures, such as the dynamic contagion process introduced in \citet{dassios2011dynamic}. The dynamic contagion process is a point process that is both self-excited and externally excited, meaning that its intensity experiences jumps driven by its own activity and by the activity of an exogenous Poisson process. In this way, the process combines the behavior of the Hawkes and shot noise processes. Hence, its infinitesimal generator can be written using a combination of expressions used in Subsections~\ref{subsecHawkes} and~\ref{subsecSN}, implying that all moments of the process can be calculated through this methodology. Similarly, these methods can also be readily applied to processes that combine dynamics from Hawkes processes and from It\^o diffusions, such as affine point processes. These processes, studied in e.g. \citet{errais2010affine,zhang2015affine,gao2019affine}, feature both self-excitement and diffusion behavior and thus have an infinitesimal generator that can be expressed using terms from the generators for Hawkes and It\^o processes. Similarly, one could study the combination of externally driven jumps and diffusive behavior such as in affine jump diffusions, see e.g. \citet{duffie2000transform}. Of course, one can also consider permutations of the model features seen in our examples, such as trading fixed size jumps for random ones to form marked Hawkes processes or changing to randomly sized batches of arrivals in the ephemerally self-exciting process. In general, the key requirement from the assumptions in Theorem~\ref{singlesol} is the closure of the system of moment differential equations specified in Equation~\ref{momentODE}. This is equivalent to having the infinitesimal generator of any polynomial being a polynomial of order no more than the original. That is, infinitesimal generators of the form
\begin{align*}
\mathcal{L}f(X_t)
&=
\underbrace{(\alpha_0 + \alpha_1 X_t ) \left(f(X_t + A_i) - f(X_t) \right)}_{\text{Up-jumps}}
+
\underbrace{(\alpha_2 + \alpha_3 X_t ) \left(f(X_t - B_i) - f(X_t) \right)}_{\text{Down-jumps}}
\\
&
\quad
+
\underbrace{(\alpha_4 + \alpha_5 X_t ) \frac{\mathrm{d}f(X_t)}{\mathrm{d}X_t}}_{\text{Drift, decay, or growth}}
+
\underbrace{(\alpha_6 + \alpha_7 X_t + \alpha_8 X_t^2 ) \frac{\mathrm{d}^2f(X_t)}{\mathrm{d}X_t^2}}_{\text{Diffusion}}
+
\underbrace{\alpha_9 \left( f(C_i X_t) - f(X_t) \right)}_{\text{Expansion or collapse}}
,
\end{align*}
can be handled by this methodology, where $\alpha_j \in \mathbb{R}$ for all $j$ and where the sequences $\{A_i\}$, $\{B_i\}$ and $\{C_i\}$ are of mutually independent random variables. Finally we note that this example generator need not be exclusive, as it is possible that other dynamics may also meet the closure requirements in Equation~\ref{momentODE}.


\section{Complexity Analysis and Numerical Experiments}\label{secNum}

\edit{
In this section we address the calculations within this method through comparison with parsimonious differential equation techniques.
Specifically, we compare both the result and the calculation time of these two methods in computing the transient moments of these stochastic processes. For the Matryoshkan-based approach, we define the calculation time as the time needed to complete the matrix computations of the moments at the specified point in time. In the differential equation approach, we take the calculation time as the time needed to reach the specified time through applying Euler's iterative method to the linear system of ODE's in Equation~\ref{genDE}, starting at time 0. We choose to compare to Euler's method because it is the fastest technique available in terms of total run time. Of course in practice in solving differential equations one may use more sophisticated approaches such as higher order Runge-Kutta methods, but such techniques are inherently at least as computationally demanding as Euler's method. As we will see however, our direct Matryoshkan matrix approach is more efficient than Euler's method outside of very short time intervals, while also delivering much more accurate answers. In Subsection~\ref{subsecComplex} we compare these two in a formal complexity analysis, and in Subsection~\ref{subsecEmp} we compare empirically through numerical experiments. Because the same matrices are used in both approaches, the time to form the matrices is omitted from each empirical calculation time reported in Subsection~\ref{subsecEmp}, although we do address the complexity of this pre-computation step in Subsection~\ref{subsecComplex}.
}

\edit{
\subsection{Complexity Characterization}\label{subsecComplex}
}

\edit{
In computing the moments through either the differential equations approach or the Matryoshkan matrix approach, one must first form the matrix that describes the system of differential equations. Thus, before comparing the two methods let us first quickly discuss this common pre-computation step. To form the matrix needed to solve for the $n^\text{th}$ moment, there will be $O(n^2)$ operations. Specifically, there are $\frac{1}{2}n(n+1)$ elements to write in this lower triangular matrix, which can be naturally conducted by writing vectors of size increasing from 1 to $n$ through the nested Matryoshkan structure of the matrix.
}

\edit{
Turning now to Euler's method, we will let $\Delta > 0$ be the time-step size parameter. For simplicity we will assume that the desired time point $t$ is a multiple of $\Delta$, meaning that there will be $t\slash \Delta$ iterations within the method to calculate the moment across time from the known initial value. Since the matrix multiplications will take $O(n^2)$ time at each step, the total complexity of Euler's method will thus be $O(n^2 t \slash \Delta)$, with the global error known to be $O(\Delta)$ \citep{butcher2016numerical}. By comparison in the direct Matryoshkan calculations, the matrix exponential calculations imply that this method is $O(n^3)$. There is of course no dependence on $\Delta$, and there exist methods in which the coefficients on $n^3$ do not depend on $t$, but rather only at lower powers of $n$ \citep{moler2003nineteen}. It is also possible this calculation could be expedited, at least in terms of the hidden coefficients, through leveraging Proposition~\ref{matryoshkanProp} and the triangularity of the matrices.
}

\edit{
In many applications, we would expect the number of Eulerian time steps should be much larger than the order of the largest moment, i.e.~$t \slash \Delta \gg n$. For example, even at a very high order moment like $n = 100$, taking a modest time step of $\Delta = .01$ would then put times $t \geq 1$ as at least as expensive for Euler's method as for the Matryoshkan approach. Of course, as the time step becomes more refined or as longer time horizons are considered, the Matryoshkan matrix calculations should become even more competitive by comparison. As we have discussed, this superiority in complexity should also immediately extend to comparisons to other numerical differential equation techniques that are themselves more complex than Euler's method. Furthermore, differential equation techniques should incur an time-step dependent error in their solution, and the closed form solutions of the Matryoshkan will not be subject to this.
}

\subsection{Empirical Comparisons and Speed Tests}\label{subsecEmp}

To demonstrate the computational efficiency and numerical precision of this method in practice, we now examine the five detailed examples covered in Subsections~\ref{subsecHawkes},~\ref{subsecSN},~\ref{subsecIto},~\ref{gcpSubsec}, and~\ref{subsecAQH}, apply our Matryoshkan matrix methodology, and compare its performance to solving the differential equations numerically through Euler's method.  To compare the results produced by the two methods, we give the absolute and relative error of Euler's method to the Matryoshkan method. That is, for $m_\mathsf{D}$ as the Eulerian moment differential equation solution and $m_\mathsf{M}$ as the Matryoshkan calculated moment, we define these errors as
$$
\text{Absolute Error}
=
|m_\mathsf{D} - m_\mathsf{M}|
\quad
\text{and}
\quad
\text{Relative Error}
=
\frac{|m_\mathsf{D} - m_\mathsf{M}|}{m_\mathsf{M}}
.
$$
All calculations are performed using simple MATLAB code on a 64-bit Windows machine with 16 GB of memory. In all four examples, we evaluate four different step sizes for Euler's method: 0.01, 0.001, 0.0001, and 0.00001. All time and error results presented in the following tables are found through averaging the results of 20 trial calculations per scenario.

We begin with the moments of the Hawkes process intensity, the first example we have discussed and the original motivation for this work. It is worth noting that to the best of our knowledge even the recognition of the matrix structure of the moments ODE system is a new contribution, as the highest order moment with explicitly stated differential equation is the fourth moment, presented without proof or solution in \citet{da2014Hawkes}. Similarly, the highest order moment with closed form solution (either transient or stationary) previously given in the literature appears to be the second. In Table~\ref{tableHawkes}, we give the errors and the time incurred for computing the first 4, 10, 20, and 100 moments. For this example we take a baseline intensity of $\lambda^* = 1$, an intensity jump size of $\alpha = 1$, a decay rate of $\beta = 2$, and we compute the moments for time $t = 10$. Additionally, we assume the initial value of the intensity is equal to the baseline. As can be quickly observed, the calculation time in the Matryoshkan computation outpaces Euler's method in all scenarios regardless of the Eulerian step size. Furthermore, when Euler's method is performed with a smaller, more precise step size, we find that its solution is increasingly close to the solution given by the Matryoshkan method as both the relative error and the absolute error decrease with the step size. In the most precise setting for Euler's method, the Matryoshkan method's run time is 4 orders of magnitude smaller for the 4, 10, and 20 moment settings and 3 orders of magnitude smaller for the $100^\text{th}$ moment.

\begin{table}[h]
\centering
\begin{tabular}{c | c c c}
 $n=4$   &  Run Time & Absolute Error & Relative Error \\
  \hline
Matryoshkan & $2.1\times10^{-4}$ sec & $\cdot$ & $\cdot$ \\
Euler $\Delta = 10^{-2}$ & $1.4\times10^{-3}$ sec & $3.0\times10^{-4}$  & $5.0\times10^{-6}$     \\
Euler $\Delta = 10^{-3}$ & $1.3\times10^{-2}$ sec & $3.1\times10^{-5}$ & $5.1\times10^{-7}$   \\
Euler $\Delta = 10^{-4}$ & $1.3\times10^{-1}$ sec & $3.1\times10^{-6}$ & $5.2\times10^{-8}$   \\
Euler $\Delta = 10^{-5}$ & $1.3\times10^{0}$ sec & $3.1\times10^{-7}$ & $5.2\times10^{-9}$   \\
\hline
\end{tabular}
\begin{tabular}{c | c c c}
$n=10$    &  Run Time & Absolute Error & Relative Error \\
  \hline
Matryoshkan  & $2.3\times10^{-4}$ sec & $\cdot$ & $\cdot$  \\
Euler $\Delta = 10^{-2}$  & $1.5\times10^{-3}$ sec & $4.4\times10^{1}$ & $1.0\times10^{-5}$  \\
Euler $\Delta = 10^{-3}$  & $1.4\times10^{-2}$ sec  & $4.5\times10^{0}$ & $1.1\times10^{-6}$  \\
Euler $\Delta = 10^{-4}$  & $1.4\times10^{-1}$ sec  & $4.5\times10^{-1}$ & $1.1\times10^{-7}$  \\
Euler $\Delta = 10^{-5}$  & $1.4\times10^{0}$ sec  & $4.5\times10^{-2}$ & $1.1\times10^{-8}$ \\
\hline
\end{tabular}
\begin{tabular}{c | c c c}
 $n = 20$   &  Run Time & Absolute Error & Relative Error \\
  \hline
Matryoshkan & $2.5\times10^{-4}$ sec & $\cdot$ & $\cdot$\\
Euler $\Delta = 10^{-2}$ &  $1.5\times10^{-3}$ sec & $7.8\times10^{12}$ & $1.7\times10^{-5}$ \\
Euler $\Delta = 10^{-3}$ &    $1.5\times10^{-2}$ sec & $8.0\times10^{11}$ & $1.8\times10^{-6}$\\
Euler $\Delta = 10^{-4}$ &   $1.4\times10^{-1}$ sec  & $8.0\times10^{10}$ & $1.8\times10^{-7}$ \\
Euler $\Delta = 10^{-5}$ &  $1.4\times10^{0}$ sec  & $8.0\times10^{9}$ & $1.8\times10^{-8}$ \\
\hline
\end{tabular}
\begin{tabular}{c | c c c}
 $n=100$  &  Run Time & Absolute Error & Relative Error \\
  \hline
Matryoshkan  & $4.3\times10^{-3}$ sec & $\cdot$ & $\cdot$\\
Euler $\Delta = 10^{-2}$ & $5.0\times10^{-3}$ sec & $3.0\times10^{145}$ & $5.1\times10^{-5}$\\
Euler $\Delta = 10^{-3}$ & $4.9\times10^{-2}$ sec & $3.1\times10^{144}$ & $5.2\times10^{-6}$\\
Euler $\Delta = 10^{-4}$ &  $4.6\times10^{-1}$ sec & $3.1\times10^{143}$ & $5.2\times10^{-7}$\\
Euler $\Delta = 10^{-5}$ &  $4.4\times10^{0}$ sec & $3.1\times10^{142}$ & $5.3\times10^{-8}$
\end{tabular}
\caption{Comparison of run time and errors for Hawkes process moment calculation via Matryoshkan matrix method and via Euler differential equation methods as the moment size increases.}\label{tableHawkes}
\end{table}

We find similarly effective performance for the shot noise process. In the example shown in Table~\ref{tableSN}, we suppose that the Poisson process arrival rate is $\lambda = 1$ and that the distribution of the shot noise is LogNormal with $\mu = 0$ and $\sigma = 1$. Furthermore, we assume that the exponential decay rate is $\beta = 4$ and we evaluate the moments at time $t = 5$. Because of the scale of these moments, we now perform these computations for $n  = 5$, 10, 15, and 20. Again we see that as the step size in Euler's method decreases the differences between the pair of results shrinks, although in this case the run times of the Matryoshkan method and Euler's method with step size 0.01 are of the same magnitude. Nevertheless, as Euler's method grows increasingly precise the Matryoshkan approach becomes more favorable; its run time is 3 orders of magnitude smaller than the most precise Euler's computational duration.

\begin{table}[h]
\centering
\begin{tabular}{c | c c c}
 $n=5$   &  Run Time & Absolute Error & Relative Error \\
  \hline
Matryoshkan & $2.1\times10^{-4}$ sec & $\cdot$ & $\cdot$ \\
Euler $\Delta = 10^{-2}$ & $1.9\times10^{-4}$ sec & $1.3\times10^{-9}$  & $3.2\times10^{-10}$     \\
Euler $\Delta = 10^{-3}$ & $2.0\times10^{-3}$ sec & $1.4\times10^{-10}$ & $3.5\times10^{-11}$   \\
Euler $\Delta = 10^{-4}$ & $1.9\times10^{-2}$ sec & $1.4\times10^{-11}$ & $3.5\times10^{-12}$   \\
Euler $\Delta = 10^{-5}$ & $1.9\times10^{-1}$ sec & $1.6\times10^{-12}$ & $3.9\times10^{-13}$   \\
\hline
\end{tabular}
\begin{tabular}{c | c c c}
$n=10$    &  Run Time & Absolute Error & Relative Error \\
  \hline
Matryoshkan  & $2.2\times10^{-4}$ sec & $\cdot$ & $\cdot$  \\
Euler $\Delta = 10^{-2}$  & $1.9\times10^{-4}$ sec & $2.6\times10^{2}$ & $5.6\times10^{-10}$  \\
Euler $\Delta = 10^{-3}$  & $1.9\times10^{-3}$ sec  & $2.9\times10^{1}$ & $6.0\times10^{-11}$  \\
Euler $\Delta = 10^{-4}$  & $1.9\times10^{-2}$ sec  & $2.9\times10^{0}$ & $6.1\times10^{-12}$  \\
Euler $\Delta = 10^{-5}$  & $1.9\times10^{-1}$ sec  & $3.2\times10^{-1}$ & $6.7\times10^{-13}$ \\
\hline
\end{tabular}
\begin{tabular}{c | c c c}
 $n = 15$   &  Run Time & Absolute Error & Relative Error \\
  \hline
Matryoshkan & $2.6\times10^{-4}$ sec & $\cdot$ & $\cdot$\\
Euler $\Delta = 10^{-2}$ &  $2.2\times10^{-4}$ sec & $1.6\times10^{39}$ & $8.4\times10^{-10}$ \\
Euler $\Delta = 10^{-3}$ &  $2.0\times10^{-3}$ sec & $1.8\times10^{38}$ & $9.1\times10^{-11}$\\
Euler $\Delta = 10^{-4}$ &  $2.0\times10^{-2}$ sec  & $1.8\times10^{37}$ & $9.2\times10^{-12}$ \\
Euler $\Delta = 10^{-5}$ &  $2.0\times10^{-1}$ sec  & $2.0\times10^{36}$ & $1.0\times10^{-12}$ \\
\hline
\end{tabular}
\begin{tabular}{c | c c c}
 $n=20$  &  Run Time & Absolute Error & Relative Error \\
  \hline
Matryoshkan  & $3.4\times10^{-4}$ sec & $\cdot$ & $\cdot$\\
Euler $\Delta = 10^{-2}$ & $2.1\times10^{-4}$ sec & $2.1\times10^{115}$ & $1.1\times10^{-9}$\\
Euler $\Delta = 10^{-3}$ & $2.1\times10^{-3}$ sec & $2.3\times10^{114}$ & $1.2\times10^{-10}$\\
Euler $\Delta = 10^{-4}$ &  $2.1\times10^{-2}$ sec & $2.3\times10^{113}$ & $1.2\times10^{-11}$\\
Euler $\Delta = 10^{-5}$ &  $2.1\times10^{-1}$ sec & $2.5\times10^{112}$ & $1.3\times10^{-12}$
\end{tabular}
\caption{Comparison of run time and errors for shot noise process moment calculation via Matryoshkan matrix method and via Euler differential equation methods as the moment size increases.}\label{tableSN}
\end{table}

In Table~\ref{tableCIR} we perform these computational experiments for the moments of a Cox-Ingersoll-Ross (CIR) process, which is an It\^o diffusion with parameter $\gamma = 1$. Moreover we assume that $\mu = 1$,  $\theta = 1$, and $\sigma = 1$, and we compute the moments at time $t = 5$ for $n = 4$, 10, 20, and 100. Like in the Hawkes process example, the Matryoshkan approach outperforms Euler's method in terms of calculation time in this example in all moment scenarios and step sizes. Moreover we again see that as the step size decreases the error between the two methods decreases while the Eulerian computation time becomes much larger than the Matryoshkan run time. Specifically for the first 4, 10, and 20 moment calculation experiments the Matryoshkan is faster by 4 orders of magnitude and in the $n = 100$ setting is is 3 orders of magnitude faster.

\begin{table}[h]
\centering
\begin{tabular}{c | c c c}
 $n=4$   &  Run Time & Absolute Error & Relative Error \\
  \hline
Matryoshkan & $2.3\times10^{-4}$ sec & $\cdot$ & $\cdot$ \\
Euler $\Delta = 10^{-2}$ & $1.4\times10^{-3}$ sec & $6.8\times10^{-3}$  & $9.3\times10^{-4}$     \\
Euler $\Delta = 10^{-3}$ & $1.3\times10^{-2}$ sec & $4.8\times10^{-4}$ & $6.6\times10^{-5}$   \\
Euler $\Delta = 10^{-4}$ & $1.3\times10^{-1}$ sec & $4.9\times10^{-5}$ & $6.6\times10^{-6}$   \\
Euler $\Delta = 10^{-5}$ & $1.3\times10^{0}$ sec & $6.8\times10^{-6}$ & $9.4\times10^{-7}$   \\
\hline
\end{tabular}
\begin{tabular}{c | c c c}
$n=10$    &  Run Time & Absolute Error & Relative Error \\
  \hline
Matryoshkan  & $2.4\times10^{-4}$ sec & $\cdot$ & $\cdot$  \\
Euler $\Delta = 10^{-2}$  & $1.4\times10^{-3}$ sec & $8.2\times10^{2}$ & $2.3\times10^{-3}$  \\
Euler $\Delta = 10^{-3}$  & $1.3\times10^{-2}$ sec  & $5.8\times10^{1}$ & $1.6\times10^{-4}$  \\
Euler $\Delta = 10^{-4}$  & $1.3\times10^{-1}$ sec  & $5.8\times10^{0}$ & $1.6\times10^{-5}$  \\
Euler $\Delta = 10^{-5}$  & $1.3\times10^{0}$ sec  & $8.3\times10^{-1}$ & $2.3\times10^{-6}$ \\
\hline
\end{tabular}
\begin{tabular}{c | c c c}
 $n = 20$   &  Run Time & Absolute Error & Relative Error \\
  \hline
Matryoshkan & $2.7\times10^{-4}$ sec & $\cdot$ & $\cdot$\\
Euler $\Delta = 10^{-2}$ &  $1.4\times10^{-3}$ sec & $1.8\times10^{11}$ & $4.3\times10^{-3}$ \\
Euler $\Delta = 10^{-3}$ &  $1.3\times10^{-2}$ sec & $1.3\times10^{10}$ & $2.9\times10^{-4}$\\
Euler $\Delta = 10^{-4}$ &  $1.3\times10^{-1}$ sec  & $1.3\times10^{9}$ & $2.9\times10^{-5}$ \\
Euler $\Delta = 10^{-5}$ &  $1.3\times10^{0}$ sec  & $1.8\times10^{8}$ & $4.3\times10^{-6}$ \\
\hline
\end{tabular}
\begin{tabular}{c | c c c}
 $n=100$  &  Run Time & Absolute Error & Relative Error \\
  \hline
Matryoshkan  & $2.0\times10^{-3}$ sec & $\cdot$ & $\cdot$\\
Euler $\Delta = 10^{-2}$ & $2.8\times10^{-3}$ sec & $4.8\times10^{127}$ & $1.3\times10^{-2}$\\
Euler $\Delta = 10^{-3}$ & $2.7\times10^{-2}$ sec & $2.1\times10^{126}$ & $5.6\times10^{-4}$\\
Euler $\Delta = 10^{-4}$ &  $2.8\times10^{-1}$ sec & $2.1\times10^{125}$ & $5.6\times10^{-5}$\\
Euler $\Delta = 10^{-5}$ &  $2.7\times10^{0}$ sec & $4.7\times10^{124}$ & $1.2\times10^{-5}$
\end{tabular}
\caption{Comparison of run time and errors for CIR process moment calculation via Matryoshkan matrix method and via Euler differential equation methods as the moment size increases.}\label{tableCIR}
\end{table}

We evaluate the Matryoshkan matrix method for the growth-collapse process in Table~\ref{tableGCP}. Like in Section 4 of \citet{boxma2006markovian}, we take $\lambda = 1$ and we also set $\mu = \frac{1}{2}$. We evaluate the moments at time $t = 8$ and the moments $n = 5$, 10, 15, and 20. In addition to observing that the Matryoshkan approach is an order of magnitude faster than any of the differential equation methods, we can also note that the relative errors are the largest we have seen in these experiments across all the Eulerian step sizes. At best, the relative error is of order $10^{-6}$, and in this case the Matryoshkan method run time is four orders of magnitude faster. As was the case for each of the other processes, we can observe that as the step size is decreased, the increased precision in Euler's method yields results closer and closer to the moments calculated by the Matryoshkan approach.

\begin{table}[h]
\centering
\begin{tabular}{c | c c c}
 $n=5$   &  Run Time & Absolute Error & Relative Error \\
  \hline
Matryoshkan & $2.2\times10^{-4}$ sec & $\cdot$ & $\cdot$ \\
Euler $\Delta = 10^{-2}$ & $1.9\times10^{-3}$ sec & $4.9\times10^{0}$  & $1.4\times10^{-3}$     \\
Euler $\Delta = 10^{-3}$ & $1.1\times10^{-2}$ sec & $7.3\times10^{-1}$ & $2.1\times10^{-4}$   \\
Euler $\Delta = 10^{-4}$ & $1.1\times10^{-1}$ sec & $4.9\times10^{-2}$ & $1.4\times10^{-5}$   \\
Euler $\Delta = 10^{-5}$ & $1.1\times10^{0}$ sec & $4.9\times10^{-3}$ & $1.4\times10^{-6}$   \\
\hline
\end{tabular}
\begin{tabular}{c | c c c}
$n=10$    &  Run Time & Absolute Error & Relative Error \\
  \hline
Matryoshkan  & $2.2\times10^{-4}$ sec & $\cdot$ & $\cdot$  \\
Euler $\Delta = 10^{-2}$  & $1.2\times10^{-3}$ sec & $1.1\times10^{6}$ & $1.7\times10^{-2}$  \\
Euler $\Delta = 10^{-3}$  & $1.2\times10^{-2}$ sec  & $1.6\times10^{5}$ & $2.6\times10^{-3}$  \\
Euler $\Delta = 10^{-4}$  & $1.2\times10^{-1}$ sec  & $1.1\times10^{4}$ & $1.7\times10^{-4}$  \\
Euler $\Delta = 10^{-5}$  & $1.2\times10^{0}$ sec  & $1.1\times10^{3}$ & $1.7\times10^{-5}$ \\
\hline
\end{tabular}
\begin{tabular}{c | c c c}
 $n = 15$   &  Run Time & Absolute Error & Relative Error \\
  \hline
Matryoshkan & $2.4\times10^{-4}$ sec & $\cdot$ & $\cdot$\\
Euler $\Delta = 10^{-2}$ &  $1.2\times10^{-3}$ sec & $9.7\times10^{10}$ & $6.3\times10^{-2}$ \\
Euler $\Delta = 10^{-3}$ &  $1.2\times10^{-2}$ sec & $1.2\times10^{10}$ & $7.9\times10^{-3}$\\
Euler $\Delta = 10^{-4}$ &  $1.2\times10^{-1}$ sec  & $9.9\times10^{8}$ & $6.4\times10^{-4}$ \\
Euler $\Delta = 10^{-5}$ &  $1.2\times10^{0}$ sec  & $9.9\times10^{7}$ & $6.4\times10^{-5}$ \\
\hline
\end{tabular}
\begin{tabular}{c | c c c}
 $n=20$  &  Run Time & Absolute Error & Relative Error \\
  \hline
Matryoshkan  & $2.6\times10^{-4}$ sec & $\cdot$ & $\cdot$\\
Euler $\Delta = 10^{-2}$ & $1.2\times10^{-3}$ sec & $5.7\times10^{15}$ & $1.3\times10^{-1}$\\
Euler $\Delta = 10^{-3}$ & $1.2\times10^{-2}$ sec & $6.9\times10^{14}$ & $1.6\times10^{-2}$\\
Euler $\Delta = 10^{-4}$ &  $1.2\times10^{-1}$ sec & $6.1\times10^{13}$ & $1.4\times10^{-3}$\\
Euler $\Delta = 10^{-5}$ &  $1.2\times10^{0}$ sec & $6.1\times10^{12}$ & $1.4\times10^{-4}$
\end{tabular}
\caption{Comparison of run time and errors for growth-collapse process moment calculation via Matryoshkan matrix method and via Euler differential equation methods as the moment size increases.}\label{tableGCP}
\end{table}

As a final table of computation comparisons, we now evaluate the moments of the Affine Queue-Hawkes process with baseline intensity $\nu^* = 1$, intensity jump size $\alpha = 2$, and duration rate $\mu = 3$. Table~\ref{tableAQH} contains the calculation times and errors for computing the first 4, 10, 20, and 100 moments of this process at time $t = 5$. Again a familiar pattern emerges, as the Matryoshkan method performance is comparable or better relative to Euler's method across all experiment scenarios. At the largest Eulerian step size the run times are of the same order but for each order of magnitude decrease in step size the method's step size becomes approximately a factor of 10 times slower than the Matryoshkan calculation. As this step size decreases the error between the two computations again decreases, implying that the Matryoshkan method also outperforms the differential equation approach in accuracy.

\begin{table}[h]
\centering
\begin{tabular}{c | c c c}
 $n=4$   &  Run Time & Absolute Error & Relative Error \\
  \hline
Matryoshkan & $2.0\times10^{-4}$ sec & $\cdot$ & $\cdot$ \\
Euler $\Delta = 10^{-2}$ & $3.3\times10^{-4}$ sec & $1.7\times10^{-1}$  & $7.8\times10^{-4}$     \\
Euler $\Delta = 10^{-3}$ & $3.2\times10^{-3}$ sec & $1.2\times10^{-2}$ & $5.6\times10^{-5}$   \\
Euler $\Delta = 10^{-4}$ & $3.2\times10^{-2}$ sec & $1.2\times10^{-3}$ & $5.6\times10^{-6}$   \\
Euler $\Delta = 10^{-5}$ & $3.1\times10^{-1}$ sec & $1.7\times10^{-4}$ & $7.9\times10^{-7}$   \\
\hline
\end{tabular}
\begin{tabular}{c | c c c}
$n=10$    &  Run Time & Absolute Error & Relative Error \\
  \hline
Matryoshkan  & $2.3\times10^{-4}$ sec & $\cdot$ & $\cdot$  \\
Euler $\Delta = 10^{-2}$  & $3.5\times10^{-4}$ sec & $8.5\times10^{6}$ & $1.9\times10^{-3}$  \\
Euler $\Delta = 10^{-3}$  & $3.3\times10^{-3}$ sec  & $6.1\times10^{5}$ & $1.3\times10^{-4}$  \\
Euler $\Delta = 10^{-4}$  & $3.3\times10^{-2}$ sec  & $6.1\times10^{4}$ & $1.3\times10^{-5}$  \\
Euler $\Delta = 10^{-5}$  & $3.3\times10^{-1}$ sec  & $8.6\times10^{3}$ & $1.9\times10^{-6}$ \\
\hline
\end{tabular}
\begin{tabular}{c | c c c}
 $n = 20$   &  Run Time & Absolute Error & Relative Error \\
  \hline
Matryoshkan & $2.6\times10^{-4}$ sec & $\cdot$ & $\cdot$\\
Euler $\Delta = 10^{-2}$ &  $3.6\times10^{-4}$ sec & $6.2\times10^{22}$ & $3.6\times10^{-3}$ \\
Euler $\Delta = 10^{-3}$ &    $3.5\times10^{-3}$ sec & $4.3\times10^{21}$ & $2.5\times10^{-4}$\\
Euler $\Delta = 10^{-4}$ &   $3.5\times10^{-2}$ sec  & $4.3\times10^{20}$ & $2.5\times10^{-5}$ \\
Euler $\Delta = 10^{-5}$ &  $3.5\times10^{-1}$ sec  & $6.2\times10^{19}$ & $3.6\times10^{-6}$ \\
\hline
\end{tabular}
\begin{tabular}{c | c c c}
 $n=100$  &  Run Time & Absolute Error & Relative Error \\
  \hline
Matryoshkan  & $4.3\times10^{-3}$ sec & $\cdot$ & $\cdot$\\
Euler $\Delta = 10^{-2}$ & $1.9\times10^{-3}$ sec & $5.3\times10^{193}$ & $1.2\times10^{-2}$\\
Euler $\Delta = 10^{-3}$ & $1.8\times10^{-2}$ sec & $2.7\times10^{192}$ & $6.3\times10^{-4}$\\
Euler $\Delta = 10^{-4}$ &  $1.8\times10^{-1}$ sec & $2.7\times10^{191}$ & $6.2\times10^{-5}$\\
Euler $\Delta = 10^{-5}$ &  $1.7\times10^{0}$ sec & $5.1\times10^{190}$ & $1.2\times10^{-5}$
\end{tabular}
\caption{Comparison of run time and errors for Affine Queue-Hawkes process moment calculation via Matryoshkan matrix method and via Euler differential equation methods as the moment size increases.}\label{tableAQH}
\end{table}

\section{Conclusion}\label{secConc}

In this work, we have defined a novel sequence of matrices called Matryoshkan matrices that stack like their Russian nesting doll namesakes. In doing so, we have found a computationally efficient manner of calculating the moments of a large class of Markov processes that satisfy a closure condition for the time derivatives of their transient moments. Furthermore, this has yielded closed form expressions for the transient and steady state moments of these process. Notably, this includes the intensity of the Hawkes process, for which finding an expression for the $n^\text{th}$ moment had been an open problem. Other examples we have discussed include It\^o diffusions from the mathematical finance literature and shot noise processes from the physics literature, which showcases the breadth of this methodology. Furthermore, our computational experiments have demonstrated the efficiency of computing at a point in time rather than through time, which is a key benefit of this method over traditional approaches for solving differential equations numerically.

We can note that there are many applications of this methodology that we have not explored in this paper and are thus opportunities for future work. For example, the vector form of the moments arising from this matrix based method naturally lends itself to use in the method of moments. Thus, Matryoshkan matrices have the potential to greatly simplify estimation for the myriad of Markov processes to which they apply. Additionally, this vector of solutions may also be of use in providing computationally tractable approximations of moment generating functions. That is, by a Taylor expansion one can approximate a moment generating function by a weighted sum of its moments. Because this paper's Matryoshkan matrix methods enable efficient calculation of higher order moments, this enables higher order approximations of the moment generating function.

\edit{
As another important direction of future work, we are also interested in extending these techniques to multivariate Markov processes. This is of practical relevance in many of the settings we have described, such as point processes driven by the Hawkes or shot noise process intensities. The challenge in this case arises in the fact that a moment's differential equation now depends on the lower product moments rather than just the lower moments, so the nesting structure is not as neatly organized. Nevertheless, addressing this generalization is an extension worth pursuing, as this would render these techniques even more applicable. 
}

\ACKNOWLEDGMENT{%
We acknowledge the generous support of the National Science Foundation (NSF) for Andrew Daw's Graduate Research Fellowship under grant DGE-1650441.
}

%
\begin{APPENDICES}

\edit{
\section{Proof of Proposition~\ref{matryoshkanProp}}
\proof{Proof.}
For clarity's sake and ease of reference, we will also enumerate the proofs of each statement.
\begin{enumerate}[$i$\textit{)}]
\item Suppose $\mathbf{X}_{n}$ and $\mathbf{Y}_{n}$ are each Matryoshkan matrices. Then, by Equation~\ref{matryoshkan}, we have that
\begin{align*}
\mathbf{X}_{n} + \mathbf{Y}_{n}
&
=
\begin{bmatrix}
\mathbf{X}_{n-1} & \mathbf{0}_{n-1 \times 1} \\
\mathbf{x}_{n} & x_{n,n}
\end{bmatrix}
+
\begin{bmatrix}
\mathbf{Y}_{n-1} & \mathbf{0}_{n-1 \times 1} \\
\mathbf{y}_{n} & y_{n,n}
\end{bmatrix}
=
\begin{bmatrix}
\mathbf{X}_{n-1} + \mathbf{Y}_{n-1} & \mathbf{0}_{n-1 \times 1} \\
\mathbf{x}_{n} + \mathbf{y}_{n} & x_{n,n} + y_{n,n}
\end{bmatrix}
,
\end{align*}
and
\begin{align*}
\mathbf{X}_{n} \mathbf{Y}_{n}
&
=
\begin{bmatrix}
\mathbf{X}_{n-1} & \mathbf{0}_{n-1 \times 1} \\
\mathbf{x}_{n} & x_{n,n}
\end{bmatrix}
\begin{bmatrix}
\mathbf{Y}_{n-1} & \mathbf{0}_{n-1 \times 1} \\
\mathbf{y}_{n} & y_{n,n}
\end{bmatrix}
=
\begin{bmatrix}
\mathbf{X}_{n-1}\mathbf{Y}_{n-1} & \mathbf{0}_{n-1 \times 1} \\
\mathbf{x}_{n}\mathbf{Y}_{n-1} + x_{n,n}\mathbf{y}_{n} & x_{n,n} y_{n,n}
\end{bmatrix}
.
\end{align*}
We can now again invoke Equation~\ref{matryoshkan} to observe that these forms satisfy this definition and thus are also Matryoshkan matrices.
\item Let $\mathbf{M}_{n} \in \mathbb{R}^{n \times n}$ be a Matryoshkan matrix with all non-zero diagonal elements $m_{i,i}$ for $i \in \{1, \dots, n\}$. By definition $\mathbf{M}_n$ is lower triangular and hence its eigenvalues are on its diagonal. Since all the eigenvalues are non-zero by assumption, $\mathbf{M}_n$ is invertible. Moreover, it is known that the inverse of a lower triangular matrix is lower triangular as well. Thus, we will now solve for lower triangular matrix $\mathbf{W}_{n} \in \mathbb{R}^{n \times n}$ such that
    $
    \mathbf{I}_n = \mathbf{M}_n \mathbf{W}_n
    $
    where $\mathbf{I}_n \in \mathbb{R}^{n \times n}$ is the identity. This can be written
    \begin{align*}
    \begin{bmatrix}
    \mathbf{I}_{n-1} & \mathbf{0}_{n-1 \times 1} \\
    \mathbf{0}_{1 \times n-1} & 1
    \end{bmatrix}
    =
    \mathbf{I}_n
    =
    \mathbf{M}_n \mathbf{W}_n
    =
    \begin{bmatrix}
    \mathbf{M}_{n-1} & \mathbf{0}_{n-1 \times 1} \\
    \mathbf{m}_{n} & m_{n,n}
    \end{bmatrix}
    \begin{bmatrix}
    \mathbf{A} & \mathbf{0}_{n-1 \times 1} \\
    \mathbf{b} & c
    \end{bmatrix}
    ,
    \end{align*}
    where $\mathbf{A} \in \mathbb{R}^{n-1\times n-1}$, $b \in \mathbb{R}^{1\times n -1}$, and $c \in \mathbb{R}$. Because $m_{i,i} \ne 0$ for all $i \in \{1,\dots, n-1\}$, we also know that $\mathbf{M}_{n-1}$ is non-singular. Thus, we can see that $\mathbf{A} = \mathbf{M}_{n-1}^{-1}$ from $\mathbf{M}_{n-1}\mathbf{A} = \mathbf{I}_{n-1}$. Likewise, $c m_{n,n} = 1$ implies $c = \frac{1}{m_{n,n}}$. Then, we have that
    $$
    \mathbf{0}_{1 \times n-1}
    =
    \mathbf{m}_n \mathbf{A} + m_{n,n} \mathbf{b}
    =
    \mathbf{m}_n \mathbf{M}_{n-1}^{-1} + m_{n,n} \mathbf{b} ,
    $$
    and so $\mathbf{b} = - \frac{1}{m_{n,n}}\mathbf{m}_n \mathbf{M}_{n-1}^{-1}$. This completes the solution for $\mathbf{W}_{n}$, and hence provides the inverse of $\mathbf{M}_n$.
    \item To begin, we will prove that
    $$
    \mathbf{M}_n^k
    =
    \begin{bmatrix}
    \mathbf{M}_{n-1}^k & \mathbf{0}_{n-1 \times 1} \\
    \mathbf{m}_{n}\sum_{j=0}^{k-1} \mathbf{M}_{n-1}^j m_{n,n}^{k-1-j} & m_{n,n}^k
    \end{bmatrix}
    $$
    for $k \in \mathbb{Z}^+$. We proceed by induction. The base case, $k=1$, holds by definition. Therefore we suppose that the hypothesis holds at $k$. Then, at $k+1$ we can observe that
    \begin{align*}
    \mathbf{M}_{n}^{k+1}
    &
    =
    \mathbf{M}_{n} \mathbf{M}_{n}^{k}
    \\
    &
    =
    \begin{bmatrix}
    \mathbf{M}_{n-1} & \mathbf{0}_{n-1 \times 1} \\
    \mathbf{m}_{n} & m_{n,n}
    \end{bmatrix}
    \begin{bmatrix}
    \mathbf{M}_{n-1}^k & \mathbf{0}_{n-1 \times 1} \\
    \mathbf{m}_{n}\sum_{j=0}^{k-1} \mathbf{M}_{n-1}^j m_{n,n}^{k-1-j} & m_{n,n}^k
    \end{bmatrix}
    \\
    &
    =
    \begin{bmatrix}
    \mathbf{M}_{n-1}^{k+1} & \mathbf{0}_{n-1 \times 1} \\
    \mathbf{m}_{n}\mathbf{M}_{n-1}^k + \mathbf{m}_{n}\sum_{j=0}^{k-1} \mathbf{M}_{n-1}^j m_{n,n}^{k-j} & m_{n,n}^{k+1}
    \end{bmatrix}
    \\
    &
    =
    \begin{bmatrix}
    \mathbf{M}_{n-1}^{k+1} & \mathbf{0}_{n-1 \times 1} \\
    \mathbf{m}_{n}\sum_{j=0}^{k} \mathbf{M}_{n-1}^j m_{n,n}^{k-j} & m_{n,n}^{k+1}
    \end{bmatrix}
    ,
    \end{align*}
    which completes the induction. We now observe further that for matrices $\mathbf{A} \in \mathbb{R}^{n \times n}$ and $\mathbf{B} \in \mathbb{R}^{n \times n}$ such that $\mathbf{A}\mathbf{B} = \mathbf{B}\mathbf{A}$ and $\mathbf{A} - \mathbf{B}$ is non-singular,
    $$
    \sum_{j=0}^{k-1} \mathbf{A}^j \mathbf{B}^{k-1-j} = \left(\mathbf{A} -\mathbf{B}\right)^{-1}\left(\mathbf{A}^k-\mathbf{B}^k\right)
    .
    $$
    This relationship can verified by multiplying the left-hand side by $\mathbf{A}-\mathbf{B}$:
    $$
    (\mathbf{A}-\mathbf{B})\sum_{j=0}^{k-1} \mathbf{A}^j \mathbf{B}^{k-1-j}
    =
    \sum_{j=0}^{k-1} \mathbf{A}^{j+1} \mathbf{B}^{k-1-j} - \sum_{j=0}^{k-1} \mathbf{A}^j \mathbf{B}^{k-j}
    =
    \mathbf{A}^k-\mathbf{B}^k
    .
    $$
    This allows us to observe that
    $$
    \mathbf{M}_n^k
    =
    \begin{bmatrix}
    \mathbf{M}_{n-1}^k & \mathbf{0}_{n-1 \times 1} \\
    \mathbf{m}_{n}(\mathbf{M}_{n-1} - m_{n,n}\mathbf{I})^{-1}\left(\mathbf{M}_{n-1}^k - m_{n,n}^k\mathbf{I}\right) & m_{n,n}^k
    \end{bmatrix}
    ,
    $$
    and thus
    \begin{align*}
    e^{\mathbf{M}_{n} t}
    &
    =
    \sum_{k=0}^\infty
    \frac{t^k \mathbf{M}_n^k}{k!}
    =
    \sum_{k=0}^\infty
    \frac{t^k }{k!}
    \begin{bmatrix}
    \mathbf{M}_{n-1}^k & \mathbf{0}_{n-1 \times 1} \\
    \mathbf{m}_{n}(\mathbf{M}_{n-1} - m_{n,n}\mathbf{I})^{-1}\left(\mathbf{M}_{n-1}^k - m_{n,n}^k\mathbf{I}\right) & m_{n,n}^k
    \end{bmatrix}
    \\
    &
    =
    \begin{bmatrix}
    e^{\mathbf{M}_{n-1} t} & \mathbf{0}_{n-1 \times 1} \\
    \mathbf{m}_{n}(\mathbf{M}_{n-1} - m_{n,n}\mathbf{I})^{-1}\left(e^{\mathbf{M}_{n-1} t} - e^{m_{n,n} t}\mathbf{I}\right) & e^{m_{n,n} t}
    \end{bmatrix}
    ,
    \end{align*}
    which completes the proof. Note that because $\mathbf{M}_{n-1}$ is triangular and because we have assumed $m_{1,1}, \dots, m_{n,n}$ are distinct, we know that $\mathbf{M}_{n-1} - m_{n,n}\mathbf{I}$ is invertible. 
    \item From the statement, we seek a matrix $\mathbf{A} \in \mathbb{R}^{n-1\times n-1}$, a row vector $\mathbf{b} \in \mathbb{R}^{1 \times n -1}$, and scalar $c \in \mathbb{R}$ such that
    $$
    \begin{bmatrix}
    \mathbf{M}_{n-1} & \mathbf{0}_{n-1 \times 1} \\
    \mathbf{m}_{n} & m_{n,n}
    \end{bmatrix}
    \begin{bmatrix}
    \mathbf{A} & \mathbf{0}_{n-1 \times 1} \\
    \mathbf{b} & c
    \end{bmatrix}
    =
    \begin{bmatrix}
    \mathbf{A} & \mathbf{0}_{n-1 \times 1} \\
    \mathbf{b} & c
    \end{bmatrix}
        \begin{bmatrix}
    \mathbf{D}_{n-1} & \mathbf{0}_{n-1 \times 1} \\
    \mathbf{0}_{1 \times n-1} & m_{n,n}
    \end{bmatrix}
    $$
    where $\mathbf{D}_{n-1} \in \mathbb{R}^{n-1 \times n -1}$ is a diagonal matrix with values $m_{1,1}, \dots , m_{n-1,n-1}$. From the triangular structure of $\mathbf{M}_n$, we know that $\mathbf{D}_n$ contains all the eigenvalues of $\mathbf{M}_n$. We will now solve the resulting sub-systems. From $\mathbf{M}_{n-1} \mathbf{A} = \mathbf{A} \mathbf{D}_{n-1}$, we take $\mathbf{A} = \mathbf{U}_{n-1}$. Substituting this forward, we see that
    $$
    \mathbf{m}_n \mathbf{U}_{n-1} + m_{n,n} \mathbf{b}
    =
    \mathbf{m}_n \mathbf{A} + m_{n,n} \mathbf{b}
    =
    \mathbf{b}\mathbf{D}_{n-1}
    $$
    and so $b =  \mathbf{m}_n \mathbf{U}_{n-1}(\mathbf{D}_{n-1} - m_{n,n} \mathbf{I})^{-1}$, where as in step (iii) we are justified in inverting $\mathbf{D}_{n-1} - m_{n,n} \mathbf{I}$ due to the fact that $m_{1,1}, \dots, m_{n,n}$ are distinct. Finally, we take $c = 1$, as any value will satisfy $c m_{n,n} = c m_{n,n}$.
\end{enumerate}
\Halmos
\endproof
}

\edit{
\section{Proof of Lemma~\ref{mmlemma}}
\proof{Proof.}
The vector solution in Equation~\ref{vecsol} is known and is thus displayed for reference. Expanding this expression in bracket-notation form, by use of Proposition~\ref{matryoshkanProp} this is
\begin{align*}
\begin{bmatrix}
\mathbf{s}_{n-1}(t)\\
s_n(t)
\end{bmatrix}
&=
\begin{bmatrix}
e^{\mathbf{M}_{n-1}t} & \mathbf{0}_{n-1 \times 1} \\
\mathbf{m}_{n}
\left(\mathbf{M}_{n-1} - m_{n,n}\mathbf{I}\right)^{-1}
\left(e^{\mathbf{M}_{n-1}t} - e^{m_{n,n} t}\mathbf{I}\right)
&
e^{m_{n,n} t}
\end{bmatrix}
\begin{bmatrix}
\mathbf{s}_{n-1}(0)\\
s_n(0)
\end{bmatrix}
\\
&
\quad
-
\begin{bmatrix}
\mathbf{M}_{n-1}^{-1} & \mathbf{0}_{n-1 \times 1} \\
-\frac{1}{m_{n,n}}\mathbf{m}_{n}\mathbf{M}_{n-1}^{-1} & \frac{1}{m_{n,n}}
\end{bmatrix}
\begin{bmatrix}
\mathbf{I} - e^{\mathbf{M}_{n-1}t} & \mathbf{0}_{n-1 \times 1} \\
-\mathbf{m}_{n}
\left(\mathbf{M}_{n-1} - m_{n,n}\mathbf{I}\right)^{-1}
\left(e^{\mathbf{M}_{n-1}t} - e^{m_{n,n} t}\mathbf{I}\right)
&
1 - e^{m_{n,n} t}
\end{bmatrix}
\begin{bmatrix}
\mathbf{c}_{n-1}\\
c_n
\end{bmatrix}
.
\end{align*}
Thus, we can find $s_n(t)$ by multiplying each left side of the equality by a unit row vector in the direction of the $n^\text{th}$ coordinate, which we denote $\mathbf{v}_n^\T$. This yields
\begin{align*}
s_n(t)
&=
\mathbf{v}_n^\T
\begin{bmatrix}
\mathbf{s}_{n-1}(t)\\
s_n(t)
\end{bmatrix}
\\&=
\begin{bmatrix}
\mathbf{m}_{n}
\left(\mathbf{M}_{n-1} - m_{n,n}\mathbf{I}\right)^{-1}
\left(e^{\mathbf{M}_{n-1}t} - e^{m_{n,n} t}\mathbf{I}\right)
&
e^{m_{n,n} t}
\end{bmatrix}
\begin{bmatrix}
\mathbf{s}_{n-1}(0)\\
s_n(0)
\end{bmatrix}
\\
&
\quad
-
\begin{bmatrix}
-\frac{1}{m_{n,n}}\mathbf{m}_{n}\mathbf{M}_{n-1}^{-1} & \frac{1}{m_{n,n}}
\end{bmatrix}
\begin{bmatrix}
\mathbf{I} - e^{\mathbf{M}_{n-1}t} & \mathbf{0}_{n-1 \times 1} \\
-\mathbf{m}_{n}
\left(\mathbf{M}_{n-1} - m_{n,n}\mathbf{I}\right)^{-1}
\left(e^{\mathbf{M}_{n-1}t} - e^{m_{n,n} t}\mathbf{I}\right)
&
1 - e^{m_{n,n} t}
\end{bmatrix}
\begin{bmatrix}
\mathbf{c}_{n-1}\\
c_n
\end{bmatrix}
\\&=
\begin{bmatrix}
\mathbf{m}_{n}
\left(\mathbf{M}_{n-1} - m_{n,n}\mathbf{I}\right)^{-1}
\left(e^{\mathbf{M}_{n-1}t} - e^{m_{n,n} t}\mathbf{I}\right)
&
e^{m_{n,n} t}
\end{bmatrix}
\begin{bmatrix}
\mathbf{s}_{n-1}(0)\\
s_n(0)
\end{bmatrix}
\\
&
\quad
-
\begin{bmatrix}
-\frac{1}{m_{n,n}}\mathbf{m}_{n}\mathbf{M}_{n-1}^{-1} & \frac{1}{m_{n,n}}
\end{bmatrix}
\begin{bmatrix}
\left(\mathbf{I} - e^{\mathbf{M}_{n-1}t}\right)\mathbf{c}_{n-1}\\
-\mathbf{m}_{n}
\left(\mathbf{M}_{n-1} - m_{n,n}\mathbf{I}\right)^{-1}
\left(e^{\mathbf{M}_{n-1}t} - e^{m_{n,n} t}\mathbf{I}\right) \mathbf{c}_{n-1}
+
c_n(1 - e^{m_{n,n} t})
\end{bmatrix}
.
\end{align*}
Then by taking these inner products, we receive
\begin{align*}
s_n(t)
&=
\mathbf{m}_{n}
\left(\mathbf{M}_{n-1} - m_{n,n}\mathbf{I}\right)^{-1}
\left(e^{\mathbf{M}_{n-1}t} - e^{m_{n,n} t}\mathbf{I}\right)
\mathbf{s}_{n-1}(0)
+
s_n(0) e^{m_{n,n} t}
+
\mathbf{m}_{n}\mathbf{M}_{n-1}^{-1}
\left(\mathbf{I} - e^{\mathbf{M}_{n-1}t}\right)
\frac{\mathbf{c}_{n-1}}{m_{n,n}}
\\
&
\quad
+
\mathbf{m}_{n}
\left(\mathbf{M}_{n-1} - m_{n,n}\mathbf{I}\right)^{-1}
\left(e^{\mathbf{M}_{n-1}t} - e^{m_{n,n} t}\mathbf{I}\right)
\frac{ \mathbf{c}_{n-1} }{m_{n,n}}
-
\frac{c_n}{m_{n,n}}\left(1 - e^{m_{n,n} t}\right)
,
\end{align*}
and this simplifies to the stated solution.
\Halmos
\endproof
}

\end{APPENDICES}
%
%


\bibliographystyle{informs2014} 
\bibliography{vanHawkes} 


\end{document}